\documentclass[12pt]{article}

\usepackage{float}
\usepackage{latexsym}
\usepackage{setspace}
\usepackage{amsmath}
\usepackage{amssymb}
\usepackage{latexsym}
\usepackage{amsmath}
\usepackage{amssymb}
\usepackage{natbib}
\usepackage{subfigure}
\setcitestyle{authoryear,open={(},close={)}}
\usepackage{hyperref}
\hypersetup{colorlinks,citecolor=blue}
\usepackage{graphics}
\usepackage{graphicx}
\usepackage{multirow}
\makeatletter
\def\ps@pprintTitle{%
  \let\@oddhead\@empty
  \let\@evenhead\@empty
  \let\@oddfoot\@empty
  \let\@evenfoot\@oddfoot
}
\makeatother
\makeatletter \oddsidemargin  -.5cm \evensidemargin -.5cm \textwidth
17.5cm \topmargin -.50in \textheight 23cm

\newcommand{\bd}{\begin{definition}}
\newcommand{\ed}{\end{definition}}
\newcommand{\br}{\begin{remark}}
\newcommand{\er}{\end{remark}}
\newcommand{\bea}{\begin{eqnarray}}
\newcommand{\eea}{\end{eqnarray}}
\newcommand{\beann}{\begin{eqnarray*}}
\newcommand{\eeann}{\end{eqnarray*}}
\newtheorem{theorem}{Theorem}[section]

\newtheorem{corollary}[theorem]{Corollary}
\newtheorem{remark}{Remark}[section]

\numberwithin{equation}{section}
\numberwithin{table}{section}
\title{Estimating location parameters of several exponential distributions with ordered restriction under Linex loss function}
\author{Shrajal Bajpai$^{a}$, Lakshmi Kanta Patra$^{a}$\footnote
	{\baselineskip=10pt
		~lkpatra@iitbhilai.ac.in; patralakshmi@gmail.com}, and Suchandan Kayal$^{b}$    \\
	 \small $^a$Department of Mathematics,
	 Indian Institute of Technology Bhilai, India, 491002\\
	\small $^b$Department of Mathematics, National Institute of Technology Rourkela,  India, 769008}
\begin{document}
\date{}
\maketitle
\begin{abstract}
	\noindent Some improved estimators of the location parameters of several exponential distributions with ordered restriction are derived and compared numerically using Monte Carlo simulations. Note that the two-parameter exponential distribution is very useful in different areas like survival analysis, reliability engineering and biomedical research, where products have a guaranteed failure-free operating time before failures begin to occur. In the present manuscript, we address the component-wise estimation of location parameters of $k~(\ge 2)$ exponential distributions under an asymmetric Linex loss function. The location parameter represents a minimum guaranteed period before failure. At first, we consider the estimation of the location parameters with ordered scale parameters. Next, we address the estimation of ordered location parameters. For this, we take three different cases into account as follows: $(i)$ scale parameters are known, $(ii)$ scale parameters are unknown but equal, $(iii)$ scale parameters are unknown and unequal. In these cases, we establish general inadmissibility results. Further, using the general result, the inadmissibility of the best affine equivariant estimator is proved. The improved estimators are written in explicit forms. Additionally, we show that the results for several important life-testing schemes namely $(i)$ Type-II censoring, $(ii)$ progressive type-II censoring and $(iii)$ record value data can be obtained using i.i.d sample.Finally, for each case, the Monte Carlo simulation technique is used to compare the performance of the proposed estimators based on their risk values. The numerical results reveal a significant improvement of the proposed estimators.
\end{abstract}

\section{Introduction} 
In the reliability theory and survival analysis, the two-parameter exponential distribution plays a vital role and is widely used as a lifetime distribution. The location parameter in two-parameter exponential distribution represents the guarantee time, and the reciprocal of scale parameter represents the hazard rate. So, estimating parameters of the two-parameter  exponential distribution is an important problem. In reliability engineering and various life testing experiments, we encounter situations where prior ordering of parameters is known. Some examples are as follows. 
\begin{enumerate}
	\item[(i)] Usually, electronic products have lower failure rates in their early age compared to their late age.
	\item [(ii)] The amount of time spent on an online shopping platform by a frequent user is less than that of  an user who shops online occasionally.
	\item [(iii)] The lifetime of  local computer hardware will be shorter than that of a branded computer hardware.
	\item [(iv)] Waiting time for an express train will be longer than that for a city metro.
	\item [(v)] Servicing time for passengers at an airport taken by a newly hired staff is higher than that
	of a professional and well-experienced staff.
\end{enumerate}


In addition to these examples, for more details on the estimation problems under order-restricted parameters for various probability distributions, one may refer to \cite{silvapulle2005constrained} and \cite{van2006restricted}. An application of the order-restricted estimation problem in meta-analysis can be found in \cite{taketomi2021meta}. In this paper, we will study the estimation of ordered location parameters for the several exponential distributions. Numerous authors have studied the estimation of ordered parameters of exponential distributions in the literature. \cite{pal1992order} have considered estimating ordered location parameters of  two exponential populations with an unknown scale parameter with respect to the squared error loss function. They proved that the unrestricted standard estimators are inadmissible by proposing improved estimators. \cite{vijayasree1993mixed} investigated the mixed estimators of ordered means of two exponential distributions. They have established that a subclass of the mixed estimators beats the usual estimator, the sample mean. \cite{misra1994estimation} have studied the component-wise estimation of ordered restricted location parameter with scale parameters known and unequal under the squared error loss function. They established that the unrestricted minimum risk estimators are inadmissible under ordered restrictions. The authors also studied the mixed estimators. Component-wise estimation of ordered parameters of the $k~(\ge2)$ exponential populations has been studied by \cite{vijayasree1995componentwise}. They have proved the inadmissibility of the usual estimators of ordered parameters with respect to the squared error loss function. \cite{misra2002smooth} have derived the smooth estimator of ordered scale parameters of two gamma distributions using \cite{brewster1974improving} methodology. They have established that this smooth estimator dominates the best affine equivariant estimator (BAEE). \cite{chang2002comparison} have considered the problem of estimating linear functions of ordered scale parameters of two gamma distributions. \cite{nematollahi2011admissibility} considered the class of mixed estimators to estimate the ordered scale parameter of gamma distributions. They have obtained a subclass of the mixed estimators that improves upon the usual estimators. 
\cite{jana2015estimation} have studied estimating the ordered scale parameter of two exponential distributions with a common guarantee time. They have obtained a uniformly minimum variance unbiased estimator (UMVUE). The authors proved that the restricted MLE (RMLE) is better than the MLE. Also, the authors derived estimators that dominate the UMVUE, MLE, and RMLE under a quadratic loss function. \cite{petropoulos2017estimation} dealt with the improved estimation of the ordered scale parameter of two multivariate Lomax distributions with unknown location parameters under a quadratic loss function. He proposed various estimators that improve upon the BAEE. \cite{patra2018estimating} considered estimating the common hazard rate of two exponential distributions with ordered location parameters under a general scale invariant loss function. They have obtained estimators that dominate the BAEE. For some recent works on estimation of ordered parameters of exponential distributions, we refer to \cite{bobotas2019estimation},
\cite{patra2021componentwise}, \cite{kayal2024estimating}, \cite{jena2024estimating}, and \cite{bajpai2025improved}. 
 
In the literature, very little attention has been given in estimating the ordered location parameters of  exponential distributions under the asymmetric Linex loss function. In this work, we study the problem of estimating ordered location parameters of $k~(\ge 2)$ exponential distributions.   Let $(X_{i1},\ldots, X_{in_{i}})$ be a random sample drawn from the $i$-th population $\mathcal{P}_{i},$ $i=1,\ldots,k~(k\geq2)$. The probability density function of the $i$-th population $\mathcal{P}_{i}$ is given by
\begin{eqnarray}
	f_{i}(x;\mu_{i},\sigma_{i})=\left\{\begin{array}{ll}
		\displaystyle\frac{1}{\sigma_{i}}~\exp{\left(-\displaystyle\frac{x-\mu_{i}}{\sigma_{i}}\right)},
		& \textrm{if $x>\mu_{i}$}\\
		0,& \textrm{otherwise,}
	\end{array} \right.
\end{eqnarray}
where $\mu_i \in \mathbb{R}$ and $\sigma_{i}>0$. This paper focuses on the problem of estimating the location parameter $\mu_i$ under the loss function, given by 
\begin{equation}\label{loss}
	L(t)=q\Big[e^{\displaystyle
		pt}-
	pt-1\Big],~p\neq0,~
	q>0,
\end{equation}
where  $t=\frac{\delta_i-\mu_i}{\sigma_i}$.  We assume $q=1.$
The loss function in (\ref{loss}) was first proposed by \cite{varian1975bayesian}. The constant $p$ governs the shape of the loss function. From \cite{zellner1986bayesian}, for $p = 1$, this loss function becomes significantly asymmetric, and in this case, the
overestimation is more serious than the underestimation. When $a < 0 $ and $t<0$, the Linex loss function increases exponentially, while for $t>0$ it grows almost linearly and vice versa. It can be observed that if  $|p|$ is close to zero, the Linex loss function becomes almost symmetric and can be approximated by a squared error loss function. For small values of $|p|$, by a Taylor series expansion, we have $e^{pt}-pt-1 \approx \frac{p^2t^2}{2}$, a squared error loss function. Further, the Linex loss function behaves asymmetric when $|p|$ is large.
The main contributions of this article are as follows:
\begin{enumerate}
	\item [(i)] We consider the component-wise estimation of the location parameter $\mu_i$ when scale parameters satisfy the order restriction $\sigma_1 \le \dots\le \sigma_k$. We  prove a general inadmissibility result for finding improved estimators. Further, it is shown that the BAEE of $\mu_i$ is inadmissible. 
	
	\item [(ii)] Next, we  discuss component-wise estimation of $\mu_i$ under the ordered restriction $\mu_1\le \dots \le \mu_k$. The cases when $(i)$ all $\sigma_i$'s are known, $(ii)$ $\sigma_i$'s are unknown and equal to $\sigma$, $(iii)$ all $\sigma_i$'s are unequal and unknown are considered. For each case, we have proved general inadmissibility results. Using this result, the inadmissibility of the BAEE is established.
	
	\item [(iii)] We conduct a simulation study to asses the performance of the proposed estimators. In this purpose, we generate $50,000$ random samples using MATLAB. The PRIs of the BAEE are tabulated to compare risk performance of the proposed estimators.
\end{enumerate}
The rest of the paper is organized as follows. Section \ref{sec2} derives estimators of $\mu_{i}$, $i= 1,\dots, k$ which improve upon the BAEE and MLE  under ordered scale parameters  followed by a simulation study. In Section \ref{sec3} we have proposed estimators under ordered location parameters considering various scenarios such as when scale parameters are known and when it is unknown with both equal and unequal cases. In Subsection \ref{sec3.1} we have considered the problem when scale parameter is known followed by a simulation study, and in Subsections \ref{sec3.2} and \ref{sec3.3} estimators are derived when scale parameters are unknown with equal and unequal cases, respectively. Each subsection contains simulation study corresponding to the derived results. Section \ref{sec4} discusses several life-testing schemes shows that the corresponding result can be obtained using methodologies based on i.i.d sample.  Finally,  in Section \ref{sec5} we conclude the work with a summary of the findings.
\section{Estimation of location parameters when scale parameters are ordered}\label{sec2}
In this section, we consider the estimation of $\mu_{i}$ when the scale parameters are ordered that is $\sigma_{1}\leq\ldots\leq\sigma_{k}.$  For estimating $\mu_i$, we consider the loss function as
\begin{eqnarray}\label{linloss1}
	L(\mu_i,\sigma_{i},\delta_{i})=\exp\left\{p\left(\frac{\delta_{i}-\mu_{i}}
	{\sigma_{i}}\right)\right\}-p\left(\frac{\delta_{i}-\mu_{i}}
	{\sigma_{i}}\right)-1,~p\neq0.
\end{eqnarray}
Define, for $i=1,\dots,k$
$$X_{i(i)}=\min\{X_{i1},\ldots, X_{in_{i}}\}~~\mbox{ and }~~T_i= \sum_{j=1}^{n_i}(X_{ij}-X_{i(1)}).$$ 
Denote $\underline{X}=(X_{1(1)},\dots,X_{k(1)})$ and $\underline{T}=(T_1,\dots,T_k)$. It can be easily seen that, for this case, $(\underline{X}, \underline{T})$ is a joint complete and sufficient statistic. Further $X_{1(1)}, \dots,X_{k(1)}, T_1,\dots, T_k$ are independently distributed with 

\begin{equation} \label{m1}
	X_{i(1)}\sim Exp(\mu_i,,\sigma_i/n_i)~~\mbox{ and }~~T_i \sim Gamma (n_i-1,\sigma_i). 
	\end{equation}  

First, we derive
the unrestricted BAEE of $\mu_{i}$ using principal of invariance.   For this purpose, we consider the affine group of transformations  as 
$$\mathcal{G}_{a_i,b_i}=\{g_{a_i,b_i}(x_{ij})= a_{i}x_{ij}+b_{i},~a_i>0,~b_i\in \mathbb{R}\}.$$  
Under this group of transformations, we have the following  
$$(X_{i(1)},T_i)\rightarrow (a_iX_{i(1)}+b_i,a_iT_i)~~\mbox{ and }~~(\mu_{i},\sigma_i) \rightarrow (a_i\mu_{i}+b_i,a_i\sigma_i).$$
Using standard arguments, the form of the affine equivariant estimator of $\mu_{i}$ is obtained as
\begin{eqnarray}\label{baee1}
	\delta_{{c_{i}}}(\underline{X},\underline{T})=X_{i(1)}+c_iT_{i},
\end{eqnarray}
where $c_i$ is an arbitrary constant. The risk of the estimators of
the form (\ref{baee1}) under the loss function (\ref{linloss1}) is given by
\begin{eqnarray} \label{rbaee}
	R(\mu_{i},\sigma_i,\delta_{{c_{i}}})=E\left[\exp\left\{
	p\left(\frac{X_{i(1)}+c_iT_{i}-\mu_{i}}{\sigma_{i}}\right)\right\}-
	p\left(\frac{X_{i(1)}+c_iT_{i}-\mu_{i}}{\sigma_{i}}\right)-1\right].
\end{eqnarray}
Differentiating (\ref{rbaee}) with respect to $c_i$, and then equating to zero we
get the minimizing choice as
\begin{eqnarray}
	c_{0i}=\frac{1}{p}\Big[1-\left(\frac{n_{i}}{n_{i}-p}\right)^{\frac{1}{n_{i}}}\Big],
\end{eqnarray}
provided $n_{i}>p.$ Thus, we have the following theorem, which provides the BAEE of $\mu_i$.
\begin{theorem}
	Under the Linex loss function (\ref{linloss1}), the best affine equivariant
	estimator of $\mu_{i}$ is
	${\delta_{c_{0i}}}(\underline{X},\underline{T})$, where
	$c_{0i}=\frac{1}{p}\Big[1-\left(\frac{n_{i}}
	{n_{i}-p}\right)^{\frac{1}{n_{i}}}\Big], ~i=1,\dots,k$. 
\end{theorem}
\subsection{Improving upon the best affine equivariant estimator}
In this section, we prove a general inadmissibility result for finding improved estimator of $\mu_i$. For this purpose, we consider a sub group of the affine group as $\mathcal{G}_1=\{g_{a,b_i}: g_{a,b_i}(x_{ij})=ax_{ij}+b_{i}, a>0,b_i\in\mathbb{R}\}$.  Under this transformation, the form of the $\mathcal{G}_1$-equivariant estimators is obtained as 
\begin{eqnarray}\label{impv1}
	\delta_{\phi_{i}}(\underline{X},\underline{T})&=&X_{i(1)}+T_{i}\phi_{i}\left(\frac{T_{1}}{T_{i}},
	\ldots,\frac{T_{i-1}}{T_{i}},\frac{T_{i+1}}{T_{i}},\ldots,\frac{T_{k}}{T_{i}}\right)\nonumber\\
	&=&X_{i(1)}+W_{i}\phi_{i}(\underline{W}),
\end{eqnarray}
where $\underline{W}=(W_{1},\ldots,W_{i-1},W_{i+1},\ldots,W_{k}),$
$W_{j}=T_{j}/T_{i},~j\neq i,~j=1,\ldots,k$ and $W_{i}=T_{i}.$ 
Further, the risk function of the estimator (\ref{impv1}) is
\begin{eqnarray}
	R(\mu_{i},\sigma_{i},{\delta_{\phi}}_{i})&=&E\left[\left(\exp\left\{{p\left(
		\frac{X_{i(1)}+W_{i}\phi_{i}(\underline{W})-\mu_{i}}{\sigma_{i}}\right)}\right\}-{p\left(
		\frac{X_{i(1)}+W_{i}\phi_{i}(\underline{W})-\mu_{i}}{\sigma_{i}}\right)}-1\right)\right]\nonumber\\
	&=&E^{\underline{W}}R_{1}(\mu_{i},\sigma_{i},{\delta}_{\phi_{i}}),
\end{eqnarray}
where $R_{1}(\mu_{i},\sigma_{i},{\delta_{\phi}}_{i})$
represents the conditional risk of ${\delta_{\phi}}_{i}$, given
$\underline{W}=\underline{w}$. Now, we have
\begin{eqnarray}
	R_{1}(\mu_{i},\sigma_{i},{\delta_{\phi}}_{i})=E\left[\left(\exp\left\{{p\left(
		\frac{X_{i(1)}+W_{i}\phi_{i}(\underline{W})-\mu_{i}}{\sigma_{i}}\right)}\right\}-{p\left(
		\frac{X_{i(1)}+W_{i}\phi_{i}(\underline{W})-\mu_{i}}{\sigma_{i}}\right)}-1\right)|\underline{W}
	=\underline{w}\right].
\end{eqnarray}
Further, $R_{1}$ is a convex function in $\phi_{i}$. Thus, the choice
of $\phi_{i}$ which minimizes $R_{1}$ can be obtained from
\begin{eqnarray}\label{impv2}
	E\left[W_{i}\exp\left\{
	p\left(\frac{X_{i(1)}+W_{i}\phi_{i}(\underline{W})}{\sigma_{i}}\right)\right\}|\underline{W}=\underline{w}\right]
	=\exp\left\{\displaystyle\frac{p\mu_{i}}{\sigma_{i}}\right\}E\left(W_{i}|\underline{W}=\underline{w}\right).
\end{eqnarray}
From \cite{vijayasree1995componentwise}, $W_i|\underline{W}=\underline{w}$ follows  $Gamma\left(n,\left(\frac{1}{\sigma_i}+\sum\limits_{j(\ne i)=1}^k\frac{w_j}{\sigma_j}\right)^{-1}\right)$, where $n=\sum\limits_{i=1}^kn_i$. After some simplification, from Equation (\ref{impv2}), we get
\begin{eqnarray}
	\phi_{i}(\underline{w},\underline{\sigma})=\frac{1}{p}
	\left(1+\sum_{j(\neq
		i)=1}^{k}w_{j}\frac{\sigma_{i}}{\sigma_{j}}\right)\times\left(1-\left(
	\frac{n_{i}}{n_{i}-p}\right) ^{\frac{1}{n+1}} \right).
\end{eqnarray}
Now we apply \cite{brewster1974improving} technique to find an improved estimator. For this we need to derive the infimum
and supremum of $\phi_{i}(\underline{w},\underline{\sigma})$ with respect to unknown parameters.   Now for $p\neq0,$ we  define 
\begin{eqnarray}
	\phi_{{i}_{*}}(\underline{w})=\inf_{\sigma_1\le\dots \le \sigma_k}\phi_{i}(\underline{w},\underline{\sigma})&=& \left\{\begin{array}{ll}
		{\frac{1}{p} \left(1-\left( \frac{n_{i}}{n_{i}-p}\right)
			^{\frac{1}{n+1}} \right) \left(1+\sum\limits_{j=2}^{k}w_j\right)},
		& \textrm{if $i=1$}\\
		{-\infty},& \textrm{if $i=2,\dots, k $}
	\end{array} \right. \nonumber
\end{eqnarray}
and
\begin{eqnarray}
{\phi_{i}}^{*}(\underline{w})=	\sup_{\sigma_1\le\dots \le \sigma_k}\phi_{i}(\underline{w},\underline{\sigma})&=& \left\{\begin{array}{ll}
		{\frac{1}{p} \left(1-\left( \frac{n_{i}}{n_{i}-p}\right)
			^{\frac{1}{n+1}} \right)},
		& \textrm{if $i=1$}\\
		{{\frac{1}{p} \left(1-\left( \frac{n_{i}}{n_{i}-p}\right)^{\frac{1}{n+1}}
				\right)}\left(1+\sum \limits_{j=1}^{i-1}w_j\right)},& \textrm{if $i=2,\dots, k.$}
	\end{array} \right. \nonumber
\end{eqnarray}
We further define the following function
\begin{eqnarray}
	\phi_{0i}(\underline{w}) = \left\{\begin{array}{ll}
		{\phi_{i}}_{*}(\underline{w}),
		& \textrm{if $\phi_{i}(\underline{w})\leq{\phi_{i}}_{*}(\underline{w})$}\\
		{\phi_{i}}(\underline{w}),& \textrm{if
			${\phi_{i}}_{*}(\underline{w})<\phi_{i}(\underline{w})<{\phi_{i}}^{*}(\underline{w})$}\\
		{\phi_{i}}^{*}(\underline{w}),& \textrm{if
			$\phi_{i}(\underline{w})\geq{\phi_{i}}^{*}(\underline{w})$}.
	\end{array} \right.
\end{eqnarray}

In the following theorem, we give sufficient condition to find an estimator that dominates the estimator $	\delta_{\phi_{i}}(\underline{X},\underline{T})$ given in (\ref{impv1}).
\begin{theorem}\label{Th1}
	The estimator $\delta_{\phi_{0i}(\underline{W})}$ has uniformly smaller risk than that of $\delta_{\phi_{i}}$ under the Linex loss function in (\ref{linloss1}) provided the condition  $P_{\underline{\theta}}(\phi_{i0}(\underline{W})\neq\phi_{i}(\underline{W}))>0$ holds. 
\end{theorem}

The following corollary provides an estimator which improves upon the BAEE as an immediate application of Theorem \ref{Th1}. 
\begin{corollary}
	Under the order restriction $\sigma_{1}\le \dots \le \sigma_{k},$ 
	\begin{itemize}
		\item [(i)] for estimating $\mu_1$, the estimator 
		\begin{eqnarray}
			\delta_{\phi_{01}}= \left\{\begin{array}{ll}
				X_{1(1)}+l_1T_{1,}
				& \textrm{if $c_{01} \leq  l_1$}\\
				X_{1(1)}+c_{01}T_{1},& \textrm{if
					$ l_1 <c_{01}< u_1$}\\
				X_{1(1)}+u_1T_{1},& \textrm{if
					$c_{01}\geq u_1$}
			\end{array} \right.
		\end{eqnarray} 
		has nowhere larger risk than the BAEE under the Linex loss function (\ref{linloss1}), where 
		$$u_1={\frac{1}{p} \left(1-\left( \frac{n_{1}}{n_{1}-p}\right)
			^{\frac{1}{n+1}} \right)}~\mbox{ and }~l_1=\frac{1}{p} \left(1-\left( \frac{n_{1}}{n_{1}-p}\right)
		^{\frac{1}{n+1}} \right) \left(1+\sum_{j=2}^{k}W_j\right).$$
		
		\item[(ii)] for estimating $\mu_i, i= 2, \dots, k$ the estimator
		\begin{eqnarray}
			\delta_{\phi_{0i}} = \left\{\begin{array}{ll}
				{X_{i(1)}+u_iT_{i}},
				& \textrm{if $c_{0i} \geq{u_{i}}$}\\
				X_{i(1)}+c_{0i}T_{i},& \textrm{otherwise}
			\end{array} \right.
		\end{eqnarray} 
		has nowhere larger risk than the BAEE under the Linex loss function in (\ref{linloss1}), where  $$u_i={{\frac{1}{p} \left(1-\left( \frac{n_{i}}{n_{i}-p}\right)^{\frac{1}{n+1}}
				\right)}\left(1+\sum_{j=1}^{i-1}W_j\right)}.$$
	\end{itemize}
\end{corollary}

\subsection{A simulation study}
 In this section, we conduct a simulation study to assess the risk performance of the proposed estimators for $k=2$. We generate $50,000$ random samples from two exponential distributions with sample sizes $n_1$ and $n_2$ to perform the simulation study. We calculate the percentage risk improvement (PRI) to assess the proposed estimators' risk performance. The PRI of an estimator $\delta$ with respect to an estimator $\delta_0$ is  defined as 
 $$PRI (\delta,\delta_0) = \frac{Risk(\delta_0)-Risk(\delta)}{Risk(\delta_0)} \times 100\%.$$
 We present the PRI of the BAEEs with respect to MLE in Table \ref{tab4}. In Tables \ref{tab1},  \ref{tab2} we report the PRI of $\delta_{\phi_{01}}$ and $\delta_{\phi_{02}}$ with respect to the BAEE, and in Tables \ref{tab3}, \ref{tab5} with respect to the MLE. To calculate the PRI values, we have taken different values of $\mu_{1},\mu_{2}$, $\sigma_{1},\sigma_{2}$, and $(n_1,n_2)$. It is observed that the risk functions of the estimators $\delta_{\phi_{01}}$ and $\delta_{\phi_{02}}$  depend on the unknown parameter only through $\frac{\sigma_{1}}{\sigma_{2}}$. So, without loss of generality we take $\mu_1 = \mu_{2} = 0$. The risk functions of the best affine equivariant estimator and the MLE are constant with respect to the unknown parameters. The following observations can be made from the tabulated values of PRI.
 \begin{itemize}
 	\item [(i)] The estimators $\delta_{\phi_{01}}$ and $\delta_{\phi_{02}}$ perform better than the BAEE and MLE. The PRI of $\delta_{\phi_{01}}$ increases then decreases as function of $\sigma_1/\sigma_2$. Furthermore, when $\sigma_{1}$ and $\sigma_{2}$ are close, $\delta_{\phi_{01}}$ and $\delta_{\phi_{02}}$ perform quite well. The PRI of the BAEE increases as $ p$  increases.
 	\item [(ii)] The PRIs of $\delta_{\phi_{01}}$ and $\delta_{\phi_{02}}$ relative to the BAEE decrease as sample size increases in contrast. The PRIs of $\delta_{\phi_{01}}$ and $\delta_{\phi_{02}}$ relative to MLE increase with sample size. PRIs of $\delta_{\phi_{01}}$ and $\delta_{\phi_{02}}$  increase as $p$ increases, moreover PRI of $\delta_{\phi_{02}}$ increases as $\frac{\sigma_{1}}{\sigma_{2}}$ increases. 

\end{itemize}
\begin{table}[H]
	\caption{PRIs of  the BAEEs $\delta_{01}$ and $\delta_{02}$ with respect to the MLE.}
	\vspace{0.2cm}
	\label{tab4}
	\centering
	\resizebox{18cm}{1.5cm}{
		\begin{tabular}{ccccccccc}
			\hline  \hline
		
			\multirow{2}{*}{$(n_1,n_2)$} & $p=-1$                      & $p=-0.5$                    & $p=0.5$                     & $p=1$                       & $p=-4$                      & $p=-2$                      & $p=2$                       & $p=4$                       \\ \cline{2-9}  
			
			& $(\delta_{01},\delta_{02})$ & $(\delta_{01},\delta_{02})$ & $(\delta_{01},\delta_{02})$ & $(\delta_{01},\delta_{02})$ & $(\delta_{01},\delta_{02})$ & $(\delta_{01},\delta_{02})$ & $(\delta_{01},\delta_{02})$ & $(\delta_{01},\delta_{02})$  \\ \hline \hline 
			 $(5,5)$                      & (37.25, 36.96)              & (38.61, 38.33)              & (41.75, 41.49)              & (43.58, 43.34)              & (31.13, 30.83)              & (34.88, 34.58)              & (47.93, 47.62)              & (60.84, 62.22)              \\ 
			$(5,7)$                      & (37.08, 40.39)              & (38.46, 41.44)              & (41.69, 43.77)              & (43.63, 45.09)              & (30.93, 35.36)              & (34.69, 38.51)              & (48.45, 48.10)              & (64.56, 56.37)              \\  
			$(7,6)$                      & (40.88, 39.26)              & (41.94, 40.45)              & (44.32, 43.18)              & (45.67, 44.76)              & (35.81, 33.70)              & (38.98, 37.14)              & (48.78, 48.50)              & (57.46, 59.80)              \\ 
			$(8,10)$                     & (41.57, 43.06)              & (42.50, 43.84)              & (44.56, 45.51)              & (45.71, 46.43)              & (36.95, 39.13)              & (39.86, 41.63)              & (48.28, 48.43)              & (55.18, 53.36)              \\ 
			$(15,12)$                    & (45.40, 44.47)              & (45.93, 45.11)              & (47.05, 46.48)              & (47.64, 47.20)              & (42.56, 41.10)              & (44.39, 43.26)              & (48.90, 48.77)              & (51.77, 52.44)              \\   \hline \hline 
			
	\end{tabular}}
\end{table}

\begin{table}[H]
		\label{tab1}
	\caption{PRIs of the estimators $\delta_{\phi_{01}}$ and $\delta_{\phi_{02}}$ with respect to the BAEE.}
	\vspace{0.2cm}
		\centering
\resizebox{16cm}{11.5cm}{
	\begin{tabular}{cccccc}
		\hline \hline
		\multirow{2}{*}{$(n_1,n_2)$} & \multirow{2}{*}{$\frac{\sigma_1}{\sigma_2}$}& $p=-1$                                    & $p=-0.5$                                  & $p=0.5$                                   & $p=1$                                    \\  \cline{3-6}
		&        & $(\delta_{\phi_{01}},\delta_{\phi_{02}})$ & $(\delta_{\phi_{01}},\delta_{\phi_{02}})$ & $(\delta_{\phi_{01}},\delta_{\phi_{02}})$ & $(\delta_{\phi_{01}},\delta_{\phi_{02}})$ \\ \hline \hline
		\multirow{7}{*}{$(5,5)$}   & 0.1    & (0.03, 0.00)                              & (0.03, 0.00)                              & (0.03, 0.00)                              & (0.04, 0.00)                              \\ 
		& 0.3    & (0.85, 0.17)                              & (0.88, 0.17)                              & (0.97, 0.19)                              & (1.03, 0.19)                              \\ 
		& 0.5    & (1.62, 0.91)                              & (1.69, 0.94)                               & (1.88, 1.01)                              & (2.01, 1.05)                              \\ 
		& 0.7    & (1.59, 2.15)                              & (1.67, 2.22)                              & (1.89, 2.38)                              & (2.04, 2.46)                              \\ 
		& 0.8    & (1.34, 2.85)                              & (1.42, 2.93)                              & (1.62, 3.13)                              & (1.76, 3.24)                              \\ 
		& 0.9    & (0.97, 3.52)                              & (1.04, 3.62)                              & (1.20, 3.86)                              & (1.33, 4.00)                              \\ 
		& 0.95   & (0.76, 3.83)                              & (0.81, 3.94)                              & (0.96, 4.21)                              & (1.07, 4.35)                           \\  \hline \hline
		\multirow{7}{*}{$(5,7)$}   & 0.1    & (0.01, 0.00)                              & (0.01, 0.00)                              & (0.01, 0.00)                              & (0.01, 0.00)                              \\ 
		& 0.3    & (0.55, 0.05)                              & (0.57, 0.05)                              & (0.63, 0.05)                              & (0.67, 0.05)                              \\ 
		& 0.5    & (1.56, 0.39)                              & (1.62, 0.40)                              & (1.78, 0.43)                              & (1.88, 0.44)                              \\ 
		& 0.7    & (2.01, 1.10)                              & (2.09, 1.13)                              & (2.30, 1.20)                              & (2.43, 1.23)                              \\ 
		& 0.8    & (1.92, 1.54)                              & (1.99, 1.57)                              & (2.18, 1.66)                              & (2.30, 1.71)                              \\ 
		& 0.9    & (1.64, 1.98)                              & (1.71, 2.03)                              & (1.86, 2.14)                              & (1.96, 2.20)                              \\ 
		& 0.95   & (1.45, 2.20)                              & (1.51, 2.25)                              & (1.64, 2.38)                              & (1.72, 2.45)                            \\ \hline \hline
		\multirow{7}{*}{$(7,6)$}   & 0.1    & (0.01, 0.00)                              & (0.01, 0.00)                              & (0.01, 0.00)                              & (0.01, 0.00)                              \\ 
		& 0.3    & (0.36, 0.12)                              & (0.37, 0.13)                              & (0.39, 0.14)                              & (0.40, 0.14)                              \\ 
		& 0.5    & (1.04, 0.86)                              & (1.07, 0.88)                              & (1.14, 0.95)                              & (1.18, 0.99)                              \\ 
		& 0.7    & (1.24, 2.19)                              & (1.27, 2.26)                              & (1.35, 2.41)                              & (1.40, 2.51)                              \\ 
		& 0.8    & (1.09, 2.95)                              & (1.12, 3.03)                              & (1.18, 3.23)                              & (1.21, 3.35)                              \\
		& 0.9    & (0.81, 3.68)                              & (0.83, 3.78)                              & (0.87, 4.02)                              & (0.89, 4.16)                              \\ 
		& 0.95   & (0.63, 4.01)                              & (0.64, 4.12)                              & (0.67, 4.38)                              & (0.68, 4.53)                              \\ \hline \hline
		\multirow{7}{*}{$(8,10)$}  & 0.1    & (0.00, 0.00)                              & (0.00, 0.00)                              & (0.00, 0.00)                              & (0.00, 0.00)                              \\ 
		& 0.3    & (0.17, 0.02)                              & (0.17, 0.02)                              & (0.18, 0.02)                              & (0.18, 0.02)                              \\ 
		& 0.5    & (0.88, 0.23)                              & (0.91, 0.24)                              & (0.95, 0.25)                              & (0.98, 0.25)                              \\ 
		& 0.7    & (1.57, 0.86)                              & (1.60, 0.87)                              & (1.68, 0.91)                              & (1.73, 0.93)                              \\ 
		& 0.8    & (1.64, 1.29)                              & (1.68, 1.31)                              & (1.76, 1.36)                              & (1.80, 1.38)                              \\ 
		& 0.9    & (1.50, 1.72)                              & (1.53, 1.75)                              & (1.60, 1.82)                              & (1.64, 1.85)                              \\ 
		& 0.95   & (1.36, 1.92)                              & (1.39, 1.96)                              & (1.45, 2.03)                              & (1.49, 2.07)                              \\ \hline \hline
		\multirow{7}{*}{$(15,12)$} & 0.1    & (0.00, 0.00)                              & (0.00, 0.00)                              & (0.00, 0.00)                              & (0.00, 0.00)                              \\ 
		& 0.3    & (0.03, 0.01)                              & (0.03, 0.01)                              & (0.03, 0.01)                              & (0.03, 0.01)                              \\ 
		& 0.5    & (0.32, 0.20)                              & (0.33, 0.20)                              & (0.33, 0.21)                              & (0.34, 0.21)                              \\ 
		& 0.7    & (0.83, 0.91)                              & (0.84, 0.93)                              & (0.86, 0.96)                              & (0.87, 0.98)                              \\ 
		& 0.8    & (0.94, 1.43)                              & (0.96, 1.45)                              & (0.98, 1.50)                              & (0.99, 1.53)                              \\ 
		& 0.9    & (0.92, 1.94)                              & (0.93, 1.97)                              & (0.95, 2.03)                              & (0.96, 2.07)                              \\ 
		& 0.95   & (0.84, 2.17)                              & (0.85, 2.20)                              & (0.87, 2.27)                              & (0.88, 2.30)                             \\ \hline \hline
	\end{tabular}
}
\end{table}

\begin{table}[H]
		\caption{PRIs of the estimators $\delta_{\phi_{01}}$ and $\delta_{\phi_{02}}$ with respect to the BAEE.}
		\label{tab2}
	\centering
\resizebox{16cm}{11.5cm}{
	\begin{tabular}{llllll}
		\hline \hline
	\multirow{2}{*}{$(n_1,n_2)$} & \multirow{2}{*}{$\frac{\sigma_1}{\sigma_2}$} & $p=-4$                                    & $p=-2$                                    & $p=2$                                     & $p=4$                                    \\ \cline{3-6}
		&        & $(\delta_{\phi_{01}},\delta_{\phi_{02}})$ & $(\delta_{\phi_{01}},\delta_{\phi_{02}})$ & $(\delta_{\phi_{01}},\delta_{\phi_{02}})$ & $(\delta_{\phi_{01}},\delta_{\phi_{02}})$ \vspace{0.3cm}\\ \hline \hline
		\multirow{7}{*}{$(5,5)$}   & 0.1    & (0.02, 0.00)                              & (0.02, 0.00)                              & (0.05, 0.00)                              & (0.05, 0.00)                              \\ 
		& 0.3    & (0.69, 0.14)                              & (0.78, 0.16)                              & (1.19, 0.21)                              & (1.19, 0.21)                              \\ 
		& 0.5    & (1.31, 0.78)                              & (1.49, 0.86)                              & (2.36, 1.10)                              & (2.36, 1.10)                              \\ 
		& 0.7    & (1.26, 1.86)                              & (1.45, 2.04)                              & (2.48, 2.62)                              & (2.48, 2.62)                              \\ 
		& 0.8    & (1.05, 2.46)                              & (1.22, 2.70)                              & (2.19, 3.45)                              & (2.19, 3.45)                              \\ 
		& 0.9    & (0.76, 3.05)                              & (0.88, 3.34)                              & (1.72, 4.26)                              & (1.72, 4.26)                              \\ 
		& 0.95   & (0.58, 3.33)                              & (0.68, 3.64)                              & (1.44, 4.64)                              & (1.44, 4.64)                             \\ \hline \hline
		\multirow{7}{*}{$(5,7)$}   & 0.1    & (0.01, 0.00)                              & (0.01, 0.00)                              & (0.01, 0.00)                              & (0.02, 0.00)                              \\ 
		& 0.3    & (0.44, 0.04)                              & (0.51, 0.04)                              & (0.75, 0.05)                              & (1.43, 0.06)                              \\ 
		& 0.5    & (1.28, 0.34)                              & (1.45, 0.37)                              & (2.16, 0.48)                              & (4.75, 0.60)                              \\  
		& 0.7    & (1.64, 0.97)                              & (1.87, 1.05)                              & (2.80, 1.32)                              & (7.14, 1.58)                              \\  
		& 0.8    & (1.58, 1.36)                              & (1.78, 1.47)                              & (2.66, 1.83)                              & (7.54, 2.15)                              \\  
		& 0.9    & (1.37, 1.76)                              & (1.54, 1.90)                              & (2.27, 2.35)                              & (7.49, 2.70)                              \\ 
		& 0.95   & (1.21, 1.95)                              & (1.36, 2.10)                              & (1.99, 2.61)                              & (7.32, 3.00)                            \\ \hline \hline
		\multirow{7}{*}{$(7,6)$}   & 0.1    & (0.01, 0.00)                              & (0.01, 0.00)                              & (0.01, 0.00)                              & (0.01, 0.00)                              \\ 
		& 0.3    & (0.32, 0.11)                              & (0.35, 0.12)                              & (0.43, 0.15)                              & (0.60, 0.19)                              \\ 
		& 0.5    & (0.88, 0.74)                              & (0.98, 0.81)                              & (1.29, 1.08)                              & (1.71, 1.40)                              \\ 
		& 0.7    & (1.06, 1.91)                              & (1.17, 2.08)                              & (1.49, 2.75)                              & (1.87, 3.58)                              \\ 
		& 0.8    & (0.93, 2.58)                              & (1.03, 2.81)                              & (1.29, 3.66)                              & (1.64, 4.76)                              \\ 
		& 0.9    & (0.69, 3.23)                              & (0.77, 3.51)                              & (0.94, 4.53)                              & (1.22, 5.80)                              \\ 
		& 0.95   & (0.54, 3.53)                              & (0.60, 3.83)                              & (0.71, 4.92)                              & (0.93, 6.26)                             \\ \hline \hline
		\multirow{7}{*}{$(8,10)$}  & 0.1    & (0.00, 0.00)                              & (0.00, 0.00)                              & (0.00, 0.00)                              & (0.00, 0.00)                              \\ 
		& 0.3    & (0.15, 0.02)                              & (0.16, 0.02)                              & (0.19, 0.02)                              & (0.24, 0.03)                              \\ 
		& 0.5    & (0.78, 0.21)                              & (0.85, 0.22)                              & (1.04, 0.27)                              & (1.22, 0.29)                              \\ 
		& 0.7    & (1.39, 0.78)                              & (1.50, 0.83)                              & (1.82, 0.96)                              & (2.06, 1.04)                              \\ 
		& 0.8    & (1.46, 1.17)                              & (1.57, 1.25)                              & (1.90, 1.44)                              & (2.15, 1.55)                              \\ 
		& 0.9    & (1.34, 1.57)                              & (1.44, 1.67)                              & (1.72, 1.92)                              & (1.95, 2.06)                              \\ 
		& 0.95   & (1.23, 1.76)                              & (1.31, 1.86)                              & (1.56, 2.14)                              & (1.73, 2.30)                            \\ \hline \hline
		\multirow{7}{*}{$(15,12)$} & 0.1    & (0.00, 0.00)                              & (0.00, 0.00)                              & (0.00, 0.00)                              & (0.00, 0.00)                              \\ 
		& 0.3    & (0.03, 0.01)                              & (0.03, 0.01)                              & (0.03, 0.01)                              & (0.03, 0.01)                              \\ 
		& 0.5    & (0.30, 0.18)                              & (0.32, 0.19)                              & (0.34, 0.22)                              & (0.36, 0.25)                              \\ 
		& 0.7    & (0.78, 0.84)                              & (0.81, 0.89)                              & (0.89, 1.02)                              & (0.95, 1.11)                              \\ 
		& 0.8    & (0.89, 1.32)                              & (0.92, 1.39)                              & (1.02, 1.58)                              & (1.08, 1.72)                              \\ 
		& 0.9    & (0.86, 1.80)                              & (0.90, 1.89)                              & (0.98, 2.14)                              & (1.05, 2.31)                              \\ 
		& 0.95   & (0.79, 2.01)                              & (0.82, 2.11)                              & (0.90, 2.38)                              & (0.97, 2.57)                         \\ \hline \hline
	\end{tabular}
}
\end{table}
\newpage
\begin{table}[H]
	\caption{PRIs of the estimators $\delta_{\phi_{01}}$ and $\delta_{\phi_{02}}$ with respect to the MLE.}
	\label{tab3}
	\centering
	\resizebox{16cm}{11.5cm}{
	\begin{tabular}{cccccc}
		\hline \hline
		\multirow{2}{*}{$(n_1,n_2)$} & \multirow{2}{*}{$\frac{\sigma_1}{\sigma_2}$} & $p=-1$                                    & $p=-0.5$                                  & $p=0.5$                                   & $p=1$                                     \\ \cline{3-6} 
		&                                              & $(\delta_{\phi_{01}},\delta_{\phi_{02}})$ & $(\delta_{\phi_{01}},\delta_{\phi_{02}})$ & $(\delta_{\phi_{01}},\delta_{\phi_{02}})$ & $(\delta_{\phi_{01}},\delta_{\phi_{02}})$ \\ \hline \hline
		\multirow{7}{*}{$(5,5)$}     & 0.1                                          & (37.27, 36.96)                            & (38.62, 38.33)                            & (41.77, 41.49)                            & (43.60, 43.34)                            \\ 
		& 0.3                                          & (37.78, 37.07)                            & (39.15, 38.43)                            & (42.31, 41.60)                            & (44.16, 43.45)                            \\ 
		& 0.5                                          & (38.27, 37.54)                            & (39.65, 38.91)                            & (42.84, 42.08)                            & (44.71, 43.93)                            \\ 
		& 0.7                                          & (38.25, 38.32)                            & (39.63, 39.70)                            & (42.85, 42.88)                            & (44.73, 44.74)                            \\ 
		& 0.8                                          & (38.09, 38.76)                            & (39.48, 40.13)                            & (42.69, 43.33)                            & (44.57, 45.18)                            \\ 
		& 0.9                                          & (37.87, 39.18)                            & (39.24, 40.56)                            & (42.45, 43.75)                            & (44.33, 45.61)                            \\ 
		& 0.95                                         & (37.73, 39.38)                            & (39.11, 40.76)                            & (42.31, 43.95)                            & (44.18, 45.81)                            \\ \hline \hline
		\multirow{7}{*}{$(5,7)$}     & 0.1                                          & (37.09, 40.39)                            & (37.09, 40.39)                            & (41.70, 43.77)                            & (43.63, 45.09)                            \\ 
		& 0.3                                          & (37.43, 40.42)                            & (37.43, 40.42)                            & (42.06, 43.80)                            & (44.00, 45.12)                            \\ 
		& 0.5                                          & (38.06, 40.63)                            & (38.06, 40.63)                            & (42.73, 44.01)                            & (44.69, 45.33)                            \\ 
		& 0.7                                          & (38.34, 41.05)                            & (38.34, 41.05)                            & (43.03, 44.45)                            & (45.00, 45.77)                            \\ 
		& 0.8                                          & (38.29, 41.31)                            & (38.29, 41.31)                            & (42.96, 44.71)                            & (44.93, 46.03)                            \\ 
		& 0.9                                          & (38.12, 41.58)                            & (38.12, 41.58)                            & (42.77, 44.98)                            & (44.73, 46.30)                            \\ 
		& 0.95                                         & (37.99, 41.70)                            & (37.99, 41.70)                            & (42.64, 45.11)                            & (44.60, 46.43)                            \\ \hline \hline
		\multirow{7}{*}{$(7,6)$}     & 0.1                                          & (40.89, 39.26)                            & (41.95, 40.45)                            & (44.33, 43.18)                            & (45.68, 44.76)                            \\ 
		& 0.3                                          & (41.10, 39.34)                            & (42.16, 40.53)                            & (44.54, 43.26)                            & (45.89, 44.84)                            \\ 
		& 0.5                                          & (41.50, 39.78)                            & (42.56, 40.98)                            & (44.96, 43.72)                            & (46.32, 45.30)                            \\ 
		& 0.7                                          & (41.62, 40.59)                            & (42.68, 41.80)                            & (45.08, 44.55)                            & (46.43, 46.14)                            \\ 
		& 0.8                                          & (41.53, 41.05)                            & (42.59, 42.26)                            & (44.98, 45.02)                            & (46.33, 46.61)                            \\ 
		& 0.9                                          & (41.36, 41.50)                            & (42.42, 42.70)                            & (44.81, 45.46)                            & (46.16, 47.06)                            \\ 
		& 0.95                                         & (41.26, 41.70)                            & (42.32, 42.91)                            & (44.70, 45.67)                            & (46.05, 47.26)                            \\ \hline \hline
		\multirow{7}{*}{$(8,10)$}    & 0.1                                          & (41.57, 43.06)                            & (42.50, 43.84)                            & (44.56, 45.51)                            & (45.71, 46.43)                            \\ 
		& 0.3                                          & (41.66, 43.07)                            & (42.60, 43.85)                            & (44.66, 45.53)                            & (45.81, 46.44)                            \\ 
		& 0.5                                          & (42.08, 43.20)                            & (43.02, 43.97)                            & (45.09, 45.65)                            & (46.24, 46.56)                            \\ 
		& 0.7                                          & (42.48, 43.55)                            & (43.42, 44.33)                            & (45.50, 46.01)                            & (46.64, 46.92)                            \\ 
		& 0.8                                          & (42.52, 43.80)                            & (43.47, 44.57)                            & (45.54, 46.25)                            & (46.68, 47.17)                            \\ 
		& 0.9                                          & (42.44, 44.05)                            & (43.38, 44.82)                            & (45.45, 46.50)                            & (46.60, 47.42)                            \\ 
		& 0.95                                         & (42.36, 44.16)                            & (43.30, 44.94)                            & (45.37, 46.62)                            & (46.51, 47.53)                            \\ \hline \hline
		\multirow{7}{*}{$(15,12)$}   & 0.1                                          & (45.40, 44.47)                            & (45.93, 45.11)                            & (47.05, 46.48)                            & (47.64, 47.20)                            \\ 
		& 0.3                                          & (45.41, 44.47)                            & (45.94, 45.12)                            & (47.07, 46.48)                            & (47.66, 47.21)                            \\ 
		& 0.5                                          & (45.57, 44.58)                            & (46.10, 45.22)                            & (47.23, 46.59)                            & (47.82, 47.32)                            \\ 
		& 0.7                                          & (45.85, 44.98)                            & (46.38, 45.62)                            & (47.50, 46.99)                            & (48.10, 47.72)                            \\ 
		& 0.8                                          & (45.91, 45.26)                            & (46.44, 45.91)                            & (47.57, 47.28)                            & (48.16, 48.01)                            \\ 
		& 0.9                                          & (45.90, 45.55)                            & (46.43, 46.19)                            & (47.55, 47.56)                            & (48.15, 48.30)                            \\ 
		& 0.95                                         & (45.86, 45.67)                            & (46.39, 46.32)                            & (47.51, 47.69)                            & (48.10, 48.42)                            \\ \hline \hline
	\end{tabular}

	}
\end{table}
\newpage

\begin{table}[H]
		\caption{PRIs of the estimators $\delta_{\phi_{01}}$ and $\delta_{\phi_{02}}$ with respect to the MLE.}
		\label{tab5}
	\centering
	\resizebox{16cm}{11.5cm}{
	\begin{tabular}{cccccc}
		\hline \hline
		\multirow{2}{*}{$(n_1,n_2)$} & \multirow{2}{*}{$\frac{\sigma_1}{\sigma_2}$} & $p=-4$                                    & $p=-2$                                    & $p= 2$                                    & $p=4$                                     \\  \cline{3-6} 
		&                                              & $(\delta_{\phi_{01}},\delta_{\phi_{02}})$ & $(\delta_{\phi_{01}},\delta_{\phi_{02}})$ & $(\delta_{\phi_{01}},\delta_{\phi_{02}})$ & $(\delta_{\phi_{01}},\delta_{\phi_{02}})$ \\ \hline \hline
		\multirow{7}{*}{$(5,5)$}     & 0.1                                          & (31.15, 30.83)                            & (34.90, 34.58)                            & (48.55, 47.73)                            & (60.87, 62.22)                            \\ 
		& 0.3                                          & (31.61, 30.93)                            & (35.39, 34.69)                            & (49.16, 48.20)                            & (61.70, 62.26)                            \\ 
		& 0.5                                          & (32.03, 31.37)                            & (35.85, 35.15)                            & (49.22, 49.00)                            & (62.62, 62.43)                            \\ 
		& 0.7                                          & (32.00, 32.11)                            & (35.83, 35.92)                            & (49.07, 49.43)                            & (63.04, 62.79)                            \\ 
		& 0.8                                          & (31.86, 32.53)                            & (35.68, 36.35)                            & (48.83, 49.86)                            & (63.08, 63.05)                            \\ 
		& 0.9                                          & (31.65, 32.94)                            & (35.46, 36.76)                            & (48.68, 50.05)                            & (62.91, 63.27)                            \\ 
		& 0.95                                         & (31.53, 33.13)                            & (35.33, 36.96)                            & (48.55, 47.73)                            & (62.71, 63.37)                            \\ \hline \hline 
		\multirow{7}{*}{$(5,7)$}     & 0.1                                          & (30.93, 35.36)                            & (34.69, 38.51)                            & (48.45, 48.10)                            & (64.56, 56.37)                            \\ 
		& 0.3                                          & (31.24, 35.39)                            & (35.02, 38.53)                            & (48.84, 48.12)                            & (65.06, 56.40)                            \\ 
		& 0.5                                          & (31.81, 35.58)                            & (35.63, 38.74)                            & (49.56, 48.34)                            & (66.24, 56.64)                            \\ 
		& 0.7                                          & (32.07, 35.99)                            & (35.91, 39.16)                            & (49.89, 48.78)                            & (67.09, 57.06)                            \\ 
		& 0.8                                          & (32.02, 36.24)                            & (35.85, 39.41)                            & (49.82, 49.05)                            & (67.23, 57.31)                            \\ 
		& 0.9                                          & (31.87, 36.50)                            & (35.69, 39.67)                            & (49.61, 49.31)                            & (67.21, 57.55)                            \\ 
		& 0.95                                         & (31.77, 36.62)                            & (35.58, 39.80)                            & (49.47, 49.45)                            & (67.15, 57.68)                            \\ \hline \hline
		\multirow{7}{*}{$(7,6)$}     & 0.1                                          & (35.82, 33.70)                            & (38.98, 37.14)                            & (48.79, 48.50)                            & (57.47, 59.80)                            \\ 
		& 0.3                                          & (36.02, 33.78)                            & (39.19, 37.22)                            & (49.01, 48.58)                            & (57.72, 59.87)                            \\ 
		& 0.5                                          & (36.38, 34.19)                            & (39.57, 37.65)                            & (49.44, 49.06)                            & (58.19, 60.36)                            \\ 
		& 0.7                                          & (36.49, 34.97)                            & (39.69, 38.45)                            & (49.55, 49.92)                            & (58.26, 61.23)                            \\ 
		& 0.8                                          & (36.41, 35.41)                            & (39.61, 38.91)                            & (49.44, 50.39)                            & (58.16, 61.71)                            \\ 
		& 0.9                                          & (36.26, 35.85)                            & (39.45, 39.35)                            & (49.26, 50.83)                            & (57.98, 62.13)                            \\ 
		& 0.95                                         & (36.16, 36.04)                            & (39.34, 39.55)                            & (49.15, 51.04)                            & (57.86, 62.31)                            \\ \hline \hline
		\multirow{7}{*}{$(8,10)$}    & 0.1                                          & (36.95, 39.13)                            & (39.86, 41.63)                            & (48.28, 48.43)                            & (55.18, 53.36)                            \\ 
		& 0.3                                          & (37.04, 39.14)                            & (39.95, 41.64)                            & (48.38, 48.44)                            & (55.29, 53.37)                            \\ 
		& 0.5                                          & (37.44, 39.26)                            & (40.37, 41.76)                            & (48.82, 48.57)                            & (55.73, 53.50)                            \\ 
		& 0.7                                          & (37.83, 39.61)                            & (40.76, 42.11)                            & (49.22, 48.93)                            & (56.10, 53.85)                            \\ 
		& 0.8                                          & (37.87, 39.85)                            & (40.80, 42.36)                            & (49.26, 49.17)                            & (56.14, 54.08)                            \\ 
		& 0.9                                          & (37.79, 40.09)                            & (40.72, 42.60)                            & (49.17, 49.42)                            & (56.06, 54.32)                            \\ 
		& 0.95                                         & (37.72, 40.20)                            & (40.65, 42.72)                            & (49.08, 49.54)                            & (55.96, 54.44)                            \\ \hline \hline
		\multirow{7}{*}{$(15,12)$}   & 0.1                                          & (42.56, 41.10)                            & (44.39, 43.26)                            & (48.90, 48.77)                            & (51.77, 52.44)                            \\ 
		& 0.3                                          & (42.57, 41.11)                            & (44.41, 43.27)                            & (48.92, 48.77)                            & (51.78, 52.45)                            \\ 
		& 0.5                                          & (42.73, 41.21)                            & (44.56, 43.37)                            & (49.08, 48.88)                            & (51.94, 52.56)                            \\ 
		& 0.7                                          & (43.00, 41.60)                            & (44.84, 43.77)                            & (49.36, 49.29)                            & (52.23, 52.97)                            \\ 
		& 0.8                                          & (43.07, 41.88)                            & (44.90, 44.05)                            & (49.42, 49.58)                            & (52.29, 53.26)                            \\ 
		& 0.9                                          & (43.05, 42.16)                            & (44.89, 44.33)                            & (49.41, 49.87)                            & (52.27, 53.54)                            \\ 
		& 0.95                                         & (43.01, 42.28)                            & (44.85, 44.46)                            & (49.36, 49.99)                            & (52.23, 53.66)                            \\ \hline \hline
	\end{tabular}}
\end{table}

\section{Estimation of ordered location parameters}\label{sec3}
In this section, we study the estimation of $\mu_{i},~i=1,\dots,k$ with the
order restriction $\mu_{1}\leq\ldots\leq\mu_{k}$. We consider three different cases, such as 
{\bf(i)} all $\sigma_i$'s are known, {\bf(ii)} all $\sigma_i$'s are unknown but equal,
and {\bf(iii)} all $\sigma_i$'s are unknown and not necessarily equal.
\subsection{Improved estimation of $\mu_i$ when $\sigma_{i}$'s are known}\label{sec3.1}
Here, we  consider the problem of finding improved estimators of $\mu_i,~i=1,\dots,k$ with respect to the loss function 
\begin{eqnarray}\label{linloos2}
	L(\mu_{i},\delta_{i})=\exp\left\{
	p\left(\frac{\delta_{i}-\mu_{i}}{\sigma_{i}}\right)\right\}
	-p\left(\frac{\delta_{i}-\mu_{i}}{\sigma_{i}}\right)-1,~p\neq0,
\end{eqnarray}
when $\sigma_1,\dots,\sigma_k$ are known.  In this case, $\underline{X}$ is a complete and sufficient statistic, and  $X_{1(1)}, \dots,X_{k(1)}$ are independently distributed. Consider the location group of transformations as  $\mathcal{H}_{a_i}=\{g_{a_i}: g_{a_{i}}(x_{ij})=x_{ij}+a_{i},~j=1,\ldots,n_{i}, i=1,\ldots,k\}$. Under this location group of transformations, the form of a location equivariant estimator is obtained as
\begin{eqnarray}
	\delta_{\alpha_i}(\underline{X})=X_{i(1)}+\alpha_i,
\end{eqnarray}
where $\alpha_i$ is a real constant. The following theorem gives the unrestricted
best location equivariant estimator of $\mu_{i}.$

\begin{theorem}
	Under the Linex loss function in $(\ref{linloos2}),$ the unrestricted best location equivariant estimator of $\mu_{i}$ is $\delta_{\alpha_{0i}}(\underline{X})$, where $\alpha_{0i}=\frac{1}{p}\ln\left(\frac{n_{i}-p\sigma_{i}}{n_{i}}\right)$, with $p<n_{i}/\sigma_{i}.$
\end{theorem}
\noindent{\bf Proof:} The proof  is straightforward, and hence it is omitted. \\

Next, we aim to find an estimator that dominates $\delta_{\alpha_{0i}}(\underline{X})$ for estimating $\mu_i$. For this purpose, consider the estimators of the form
\begin{eqnarray}\label{deltapsi}
	\delta_{\psi_{i}}(\underline{X})=X_{i(1)}+\psi_{i}(\underline{Y}),
\end{eqnarray}
where
$\underline{Y}=(Y_{1},\ldots,Y_{i-1},Y_{i+1},\ldots,Y_{k}),~Y_{i}=X_{i(1)}$
and $Y_{j}=X_{j(1)}-X_{i(1)},~j\neq i,~j=1,\ldots,k$ (see \cite{vijayasree1995componentwise}).  The risk function of
the estimators $\delta_{\psi_{i}}(\underline{X})$ is
\begin{eqnarray*}
	R(\mu_{i},\delta_{\psi_{i}}(\underline{X}))&=&E\left[\exp\left\{\displaystyle
	p\left(\frac{X_{i(1)}+\psi_{i}(\underline{Y})-\mu_{i}}{\sigma_{i}}\right)\right\}-
	p\left(\frac{X_{i(1)}+\psi_{i}(\underline{Y})-\mu_{i}}{\sigma_{i}}\right)-1\right]\\
	&=&E^{\underline{Y}}R_{1}(\mu_{i},\delta_{\psi_{i}}(\underline{X})),
\end{eqnarray*}
where
$$R_{1}(\mu_{i},\delta_{\psi_{i}}(\underline{X}))=E\Big[\left(\exp\left\{
p\left(\frac{X_{i(1)}+\psi_{i}(\underline{Y})-\mu_{i}}{\sigma_{i}}\right)\right\}-
p\left(\displaystyle\frac{X_{i(1)}+\psi_{i}(\underline{Y})-\mu_{i}}{\sigma_{i}}\right)-1\right)|\underline{Y}
=\underline{y}\Big].$$ 
After some calculations, the conditional risk $R_{1}$ is minimized at
\begin{eqnarray}\label{eq3.4}
	\psi_{i}(\underline{y},\underline{\mu})=
	\frac{\sigma_{i}}{p}\ln\left[\frac{\exp\{\frac{p}{\sigma_{i}}\mu_{i}\}}
	{E\left(\exp\left\{\frac{pX_{i(1)}}{\sigma_{i}}\right\}|\underline{Y}=\underline{y}\right)}\right].
\end{eqnarray}
The conditional distribution of $X_{i(1)}$ given
$\underline{Y}=\underline{y}$ is exponential with
location parameter $\mu_0$ and scale
parameter $q^{-1},$ where
$$\mu_0=\max\left\{\mu_{i},\max_{1\leq j\leq
	k,j\neq i}(\mu_{j}-y_{j})\right\} ~~\mbox{and }~~q=\sum_{j=1}^{k}n_{j}\sigma_{j}^{-1}.$$ 
Simplifying $(\ref{eq3.4}),$ we get
\begin{eqnarray*}
	\psi_{i}(\underline{y},\underline{\mu})=
	\mu_{i}-\mu_0
	+\frac{\sigma_{i}}{p}\ln\left(\frac{\sigma_{i}q-p}{q\sigma_{i}}\right),
\end{eqnarray*}
provided $q\sigma_i>p$.
Under the order restriction $\mu_{1}\leq\ldots\leq\mu_{k},$ the
infimum and supremum of $\psi_{i}(\underline{y},\underline{\mu})$ are obtained as
\begin{eqnarray*}
	\inf_{\underline{\mu}}\psi_{i}(\underline{y},\underline{\mu}) &=&
	\left\{\begin{array}{ll} -\infty,
		& \textrm{if $i=1,\ldots,k-1$}\\
		\displaystyle\frac{\sigma_{k}}{p}\ln\left(\frac{q\sigma_{k}-p}{q\sigma_{k}}\right)
		+\min(0,y_{1},\ldots,y_{k-1}),&
		\textrm{if $i=k$}
	\end{array} \right.\nonumber\\
	&=&{\psi}_{i*}(\underline{y}),~\mbox{(say)}
\end{eqnarray*}
and
\begin{eqnarray*}
	\sup_{\underline{\mu}}\psi_i(\underline{y},\underline{\mu}) &=&
	\left\{\begin{array}{ll}
		\displaystyle\frac{\sigma_{i}}{p}\ln\left(\frac{q\sigma_{i}-p}{q\sigma_{i}}\right)+\min(0,y_{i+1},\ldots,y_{k}),
		& \textrm{if $i=1,\ldots,k-1$}\\
		\displaystyle\frac{\sigma_{k}}{p}\ln\left(\frac{q\sigma_{k}-p}{q\sigma_{k}}\right),&
		\textrm{if $i=k$}
	\end{array} \right.\nonumber\\
	&=&{\psi}_i^{*}(\underline{y}),~\mbox{(say).}
\end{eqnarray*}
Deefine
\begin{eqnarray}\label{equ22}
	\psi_{0i}(\underline{y}) = \left\{\begin{array}{ll}
		{\psi_{i}}_{*}(\underline{y}),
		& \textrm{if $\psi_{i}(\underline{y})\leq{\psi_{i}}_{*}(\underline{y})$}\\
		{\psi_{i}}(\underline{y}),& \textrm{if
			${\psi_{i}}_{*}(\underline{y})<\psi_{i}(\underline{y})<{\psi_{i}}^{*}(\underline{y})$}\\
		{\psi_{i}}^{*}(\underline{y}),& \textrm{if
			$\psi_{i}(\underline{y})\geq{\psi_{i}}^{*}(\underline{y})$}.
	\end{array} \right.
\end{eqnarray}

In the following theorem, we obtain an estimator that improves upon $\delta_{\psi_{i}}(\underline{Y})$. 
\begin{theorem}\label{theosigmak}
	Let $\delta_{\psi_{i}}(\underline{X})$ be an
	estimator of the form (\ref{deltapsi}) and $\psi_{i0}(\underline{X})$ be
	defined in (\ref{equ22}). Then, for estimating $\mu_i, i=1,\dots,k$ the estimator
	$\delta_{\psi_{0i}}(\underline{X})$ dominates $\delta_{\psi_{i}}(\underline{X})$ under the Linex loss function (\ref{linloos2}), provided 
	$P_{\underline{\mu}}(\psi_{i0}(\underline{Y})\neq\psi_{i}(\underline{Y}))>0$.
\end{theorem}

As an application of the above theorem, in the following lemma, we obtain an estimator that dominates $\delta_{\alpha_{0i}}(\underline{X})$. Denote 
$$\gamma_i=\frac{\sigma_{i}}{p}\ln\left(\frac{q\sigma_{i}-p}{q\sigma_{i}}\right)
+\min(0,Y_{i+1},\ldots,Y_{k}),~~
i=1,\ldots,k-1$$
$$\gamma_k=\frac{\sigma_{k}}{p}\ln\left(\frac{q\sigma_{k}-p}{q\sigma_{k}}\right)$$
$$\nu_k=\frac{\sigma_{k}}{p}\ln\left(\frac{q\sigma_{k}-p}{q\sigma_{k}}\right)
+\min(0,Y_{1},\ldots,Y_{k-1}).$$
\begin{corollary}
	Let $\mu_1,\dots,\mu_k$ satisfy the order restriction $\mu_{1} \le \dots \le\mu_{k}$ and $\sigma_1,\dots,\sigma_k$ are known. Then
	\begin{itemize}
		\item [(i)] for estimating $\mu_i, i=1,\dots,k-1$ the estimator 
		\begin{eqnarray}
			\delta_{\psi_{0i}}(\underline{X}) = \left\{\begin{array}{ll}
				X_{i(1)}+\gamma_i,~~\mbox{ if }~\alpha_{0i}>\gamma_i\\
				X_{i(1)}+\alpha_{0i},~~\mbox{ otherwise}
			\end{array} \right.
		\end{eqnarray}
		dominates $\delta_{\alpha_{0i}}$ under the Linex loss function (\ref{linloos2}). 
		\item[(ii)] for estimating $\mu_k$ the estimator 
		\begin{eqnarray}
			\displaystyle\delta_{\psi_{0k}}(\underline{X}) = \left\{\begin{array}{ll}
				X_{i(1)}+\nu_k,
				& \textrm{if $\alpha_{0k}\leq  \nu_k$}\\
				X_{i(1)}+\alpha_{0k}, & \textrm{if
					$\nu_k<\alpha_{0k}<\gamma_k$}\\
				X_{i(1)}+\gamma_k, & \textrm{if
					$\alpha_{0i}\geq \gamma_k$}
			\end{array} \right.
		\end{eqnarray}
		dominates the best location equivariant estimator $\delta_{\alpha_{0k}}$ under the Linex loss function (\ref{linloos2}).
	\end{itemize}
\end{corollary}
The MLE of $\mu_i$ under the order restriction $\mu_{1} \le \dots \le\mu_{k}$, that is the restricted MLE is obtained as
\begin{eqnarray}\label{RML1}\nonumber
	\delta_{iRML}&=&\min\left\{X_{i(1)},\dots,X_{k(1)}\right\}\\
	&=&X_{i(1)}+\min\left\{0,Y_{i+1},\dots,Y_k\right\},
\end{eqnarray}
for $i=1,\dots,k$.
As an immediate application of Theorem \ref{theosigmak}, in the following corollary, we prove that $\delta_{iRML}$ is inadmissible. 
\begin{corollary}
	Under the order restriction $\mu_{1} \le \dots \le\mu_{k}$, for estimating $\mu_i$  the restricted MLE $\delta_{iRML}$ is inadmissible with respect to the Linex loss function in (\ref{linloos2}) and dominated by the estimator 
	\begin{equation*}
		\delta_{iRMLI}=\min\left\{X_{i(1)},\dots,X_{k(1)}\right\}+\frac{\sigma_{i}}{p}\ln\left(\frac{q\sigma_{i}-p}{q\sigma_{i}}\right),
	\end{equation*} 
	for $i=1,\dots,k$.
\end{corollary}

\subsection{A simulation study}
In this section, we have conducted a simulation study to compare the risk performance of the proposed estimators. In Tables \ref{tab6} and \ref{tab7}, we tabulate the PRIs of $\delta_{\psi_{01}}$ and $\delta_{\psi_{02}}$ with respect to the BAEE. Tables \ref{tab1RML} and \ref{tab2RML}  provide the  PRIs of the estimator $\delta_{1RMLI}$ with respect to the restricted MLE $\delta_{1RML}$.    Tables \ref{tab10.1} and  \ref{tab10} provide the  PRIs of the estimator $\delta_{2RMLI}$ with respect to the restricted MLE $\delta_{2RML}$.
In this case, $\sigma_{1}$ and $\sigma_{2}$ are known, and by fixing the values of $\sigma_1$ and $\sigma_2$, the corresponding PRI values are calculated. The risk functions of the estimators $\delta_{\psi_{01}}$, $\delta_{\psi_{02}}$, $\delta_{1RMLI}$ depend on the unknown parameter only through $\mu_{2}-\mu_{1}$ and the risk function of the estimator $\delta_{2RMLI}$ is independent of the unknown parameter.  We have the following observations from the tabulated values. 
\begin{itemize}
	\item [(i)]  For smaller sample size, $\delta_{\psi_{01}}$ and $\delta_{\psi_{02}}$ always outperform the BAEE. When the sample size increases, the PRIs of $\delta_{\psi_{01}}$ and $\delta_{\psi_{02}}$ decrease and become very close to PRI of BAEE.
	\item [(ii)] For fixed $(n_1,n_2)$ and $\mu_{2}-\mu_{1}$, PRI values of $\delta_{\psi_{01}}$ and $\delta_{\psi_{02}}$ for  $(\sigma_{1},\sigma_{2})=(3,2)$ are higher than that of $(1,1.5)$.  This shows that the PRI values of $\delta_{\psi_{01}}$ and $\delta_{\psi_{02}}$  depend on the values of $\sigma_1$ and $\sigma_2$. 
	As the value $\mu_{2}-\mu_{1}$ increases the PRI values decrease. 
	\item[(iii)] PRIs of $\delta_{1RMLI}$ and $\delta_{2RMLI}$  increase as $p$ increases. PRI of $\delta_{2RMLI}$ decreases as sample size increases. It is observed that  $\delta_{iRMLI}$ improves upon $\delta_{iRMLI}$ significantly. 
\end{itemize}

\begin{table}[H]
	\caption{PRIs of the estimators $\delta_{\psi_{01}}$ and $\delta_{\psi_{02}}$ with respect to the BAEE when $\sigma_i$'s are known. }
	\vspace{0.2cm}
	\label{tab6}
	\centering
	\resizebox{18cm}{10cm}{
	\begin{tabular}{llllllllll}
		\hline \hline
		\multirow{3}{*}{$(n_1,n_2)$} & \multirow{3}{*}{$\mu_2-\mu_1$} & \multicolumn{4}{c}{$(\sigma_1,\sigma_2)=(1,1.5)$}                                                                               & \multicolumn{4}{c}{$(\sigma_1,\sigma_2)=(3,2)$}                                                                               \\ \cline{3-10} 
		&                              & \multicolumn{1}{l}{$p=-1$}         & \multicolumn{1}{l}{$p=-0.5$}       & \multicolumn{1}{l}{$p=0.5$}        & $p=1$          & \multicolumn{1}{l}{$p=-1$}        & \multicolumn{1}{l}{$p=-0.5$}      & \multicolumn{1}{l}{$p=0.5$}        & $p=1$          \\ \cline{3-10} 
		&                              & \multicolumn{1}{l}{{$(\delta_{\psi_{01}},\delta_{\psi_{02}})$}}              & \multicolumn{1}{l}{{$(\delta_{\psi_{01}},\delta_{\psi_{02}})$}}              & \multicolumn{1}{l}{{$(\delta_{\psi_{01}},\delta_{\psi_{02}})$}}              & {$(\delta_{\psi_{01}},\delta_{\psi_{02}})$}              & \multicolumn{1}{l}{$(\delta_{\psi_{01}},\delta_{\psi_{02}})$}             & \multicolumn{1}{l}{$(\delta_{\psi_{01}},\delta_{\psi_{02}})$}             & \multicolumn{1}{l}{{$(\delta_{\psi_{01}},\delta_{\psi_{02}})$}}              & {$(\delta_{\psi_{01}},\delta_{\psi_{02}})$}             \\ \hline \hline
		\multirow{7}{*}{$(5,5)$}    & 0.1                          & \multicolumn{1}{l}{(47.84, 11.02)} & \multicolumn{1}{l}{(48.37, 12.90)} & \multicolumn{1}{l}{(49.41, 18.42)} & (49.89, 22.63) & \multicolumn{1}{l}{(73.40, 4.73)} & \multicolumn{1}{l}{(70.37, 6.56)} & \multicolumn{1}{l}{(55.61, 13.63)} & (37.38, 20.69) \\ 
		& 0.2                          & \multicolumn{1}{l}{(37.42, 6.70)}  & \multicolumn{1}{l}{(38.42, 7.98)}  & \multicolumn{1}{l}{(40.58, 12.11)} & (41.72, 15.66) & \multicolumn{1}{l}{(70.58, 3.99)} & \multicolumn{1}{l}{(67.86, 5.55)} & \multicolumn{1}{l}{(53.98, 11.73)} & (36.39, 18.28) \\ 
		& 0.4                          & \multicolumn{1}{l}{(19.28, 2.44)}  & \multicolumn{1}{l}{(20.54, 2.89)}  & \multicolumn{1}{l}{(23.58, 4.38)}  & (25.40, 5.69)  & \multicolumn{1}{l}{(62.17, 2.86)} & \multicolumn{1}{l}{(60.30, 3.96)} & \multicolumn{1}{l}{(48.88, 8.37)}  & (33.23, 13.08) \\ 
		& 0.6                          & \multicolumn{1}{l}{(8.76, 0.91)}   & \multicolumn{1}{l}{(9.69, 1.08)}   & \multicolumn{1}{l}{(12.21, 1.64)}  & (13.89, 2.11)  & \multicolumn{1}{l}{(52.42, 2.04)} & \multicolumn{1}{l}{(51.40, 2.82)} & \multicolumn{1}{l}{(42.56, 5.94)}  & (29.31, 9.30)  \\ 
		& 0.9                          & \multicolumn{1}{l}{(2.42, 0.18)}   & \multicolumn{1}{l}{(2.83, 0.22)}   & \multicolumn{1}{l}{(4.08, 0.35)}   & (5.05, 0.46)   & \multicolumn{1}{l}{(38.35, 1.23)} & \multicolumn{1}{l}{(38.27, 1.71)} & \multicolumn{1}{l}{(33.02, 3.64)}  & (23.27, 5.69)  \\ 
		& 1.2                          & \multicolumn{1}{l}{(0.55, 0.05)}   & \multicolumn{1}{l}{(0.67, 0.06)}   & \multicolumn{1}{l}{(1.14, 0.09)}   & (1.58, 0.11)   & \multicolumn{1}{l}{(26.92, 0.74)} & \multicolumn{1}{l}{(27.38, 1.03)} & \multicolumn{1}{l}{(24.68, 2.21)}  & (17.78, 3.46)  \\ 
		& 1.4                          & \multicolumn{1}{l}{(0.22, 0.02)}   & \multicolumn{1}{l}{(0.28, 0.02)}   & \multicolumn{1}{l}{(0.55, 0.03)}   & (0.83, 0.04)   & \multicolumn{1}{l}{(20.93, 0.52)} & \multicolumn{1}{l}{(21.55, 0.73)} & \multicolumn{1}{l}{(19.96, 1.57)}  & (14.66, 2.47)  \\ \hline \hline
		\multirow{7}{*}{$(5,7)$}    & 0.1                          & \multicolumn{1}{l}{(54.00, 8.65)}  & \multicolumn{1}{l}{(54.40, 9.85)}  & \multicolumn{1}{l}{(55.07, 13.16)} & (55.07, 13.16) & \multicolumn{1}{l}{(77.40, 3.57)} & \multicolumn{1}{l}{(74.05, 4.76)} & \multicolumn{1}{l}{(58.34, 8.78)}  & (39.20, 12.26) \\ 
		& 0.2                          & \multicolumn{1}{l}{(42.56, 5.21)}  & \multicolumn{1}{l}{(43.54, 5.95)}  & \multicolumn{1}{l}{(45.58, 8.01)}  & (45.58, 8.01)  & \multicolumn{1}{l}{(74.46, 3.02)} & \multicolumn{1}{l}{(71.45, 4.01)} & \multicolumn{1}{l}{(56.65, 7.38)}  & (38.17, 10.33) \\ 
		& 0.4                          & \multicolumn{1}{l}{(22.10, 1.94)}  & \multicolumn{1}{l}{(23.45, 2.21)}  & \multicolumn{1}{l}{(26.62, 2.94)}  & (26.62, 2.94)  & \multicolumn{1}{l}{(65.63, 2.21)} & \multicolumn{1}{l}{(63.53, 2.92)} & \multicolumn{1}{l}{(51.30, 5.35)}  & (34.94, 7.47)  \\ 
		& 0.6                          & \multicolumn{1}{l}{(10.12, 0.69)}  & \multicolumn{1}{l}{(11.14, 0.79)}  & \multicolumn{1}{l}{(13.84, 1.07)}  & (13.84, 1.07)  & \multicolumn{1}{l}{(55.34, 1.52)} & \multicolumn{1}{l}{(54.15, 2.02)} & \multicolumn{1}{l}{(44.81, 3.73)}  & (30.92, 5.23)  \\ 
		& 0.9                          & \multicolumn{1}{l}{(2.76, 0.14)}   & \multicolumn{1}{l}{(3.21, 0.16)}   & \multicolumn{1}{l}{(4.53, 0.22)}   & (4.53, 0.22)   & \multicolumn{1}{l}{(40.72, 0.90)} & \multicolumn{1}{l}{(40.58, 1.21)} & \multicolumn{1}{l}{(34.98, 2.24)}  & (24.65, 3.14)  \\ 
		& 1.2                          & \multicolumn{1}{l}{(0.67, 0.04)}   & \multicolumn{1}{l}{(0.81, 0.04)}   & \multicolumn{1}{l}{(1.28, 0.05)}   & (1.28, 0.05)   & \multicolumn{1}{l}{(28.78, 0.61)} & \multicolumn{1}{l}{(29.22, 0.80)} & \multicolumn{1}{l}{(26.26, 1.44)}  & (18.93, 1.99)  \\ 
		& 1.4                          & \multicolumn{1}{l}{(0.21, 0.01)}   & \multicolumn{1}{l}{(0.26, 0.02)}   & \multicolumn{1}{l}{(0.44, 0.02)}   & (0.44, 0.02)   & \multicolumn{1}{l}{(22.44, 0.41)} & \multicolumn{1}{l}{(23.07, 0.55)} & \multicolumn{1}{l}{(21.33, 1.01)}  & (15.63, 1.42)  \\ \hline \hline
		\multirow{7}{*}{$(7,6)$}    & 0.1                          & \multicolumn{1}{l}{(41.41, 11.43)} & \multicolumn{1}{l}{(41.93, 13.14)} & \multicolumn{1}{l}{(42.99, 17.90)} & (43.52, 21.30) & \multicolumn{1}{l}{(68.85, 6.21)} & \multicolumn{1}{l}{(66.34, 8.05)} & \multicolumn{1}{l}{(56.62, 14.30)} & (56.62, 14.30) \\ 
		& 0.2                          & \multicolumn{1}{l}{(27.88, 5.67)}  & \multicolumn{1}{l}{(28.72, 6.59)}  & \multicolumn{1}{l}{(30.53, 9.44)}  & (31.51, 11.75) & \multicolumn{1}{l}{(64.32, 4.85)} & \multicolumn{1}{l}{(62.19, 6.30)} & \multicolumn{1}{l}{(53.50, 11.39)} & (53.50, 11.39) \\ 
		& 0.4                          & \multicolumn{1}{l}{(9.82, 1.40)}   & \multicolumn{1}{l}{(10.52, 1.63)}  & \multicolumn{1}{l}{(12.19, 2.32)}  & (13.19, 2.88)  & \multicolumn{1}{l}{(52.11, 3.03)} & \multicolumn{1}{l}{(50.92, 3.93)} & \multicolumn{1}{l}{(44.77, 7.07)}  & (44.77, 7.07)  \\ 
		& 0.6                          & \multicolumn{1}{l}{(2.99, 0.36)}   & \multicolumn{1}{l}{(3.33, 0.41)}   & \multicolumn{1}{l}{(4.26, 0.58)}   & (4.87, 0.72)   & \multicolumn{1}{l}{(39.72, 1.87)} & \multicolumn{1}{l}{(39.28, 2.43)} & \multicolumn{1}{l}{(35.44, 4.40)}  & (35.44, 4.40)  \\ 
		& 0.9                          & \multicolumn{1}{l}{(0.49, 0.05)}   & \multicolumn{1}{l}{(0.57, 0.06)}   & \multicolumn{1}{l}{(0.83, 0.08)}   & (1.03, 0.09)   & \multicolumn{1}{l}{(24.45, 0.94)} & \multicolumn{1}{l}{(24.62, 1.22)} & \multicolumn{1}{l}{(23.09, 2.22)}  & (23.09, 2.22)  \\ 
		& 1.2                          & \multicolumn{1}{l}{(0.06, 0.01)}   & \multicolumn{1}{l}{(0.08, 0.01)}   & \multicolumn{1}{l}{(0.12, 0.01)}   & (0.16, 0.01)   & \multicolumn{1}{l}{(14.18, 0.46)} & \multicolumn{1}{l}{(14.56, 0.60)} & \multicolumn{1}{l}{(14.24, 1.09)}  & (14.24, 1.09)  \\ 
		& 1.4                          & \multicolumn{1}{l}{(0.00, 0.00)}   & \multicolumn{1}{l}{(0.00, 0.00)}   & \multicolumn{1}{l}{(0.00, 0.00)}   & (0.00, 0.00)   & \multicolumn{1}{l}{(9.57, 0.26)}  & \multicolumn{1}{l}{(9.95, 0.35)}  & \multicolumn{1}{l}{(10.07, 0.65)}  & (10.07, 0.65)  \\ \hline \hline
		\multirow{7}{*}{$(8,10)$}   & 0.1                          & \multicolumn{1}{l}{(46.56, 7.59)}  & \multicolumn{1}{l}{(46.56, 7.59)}  & \multicolumn{1}{l}{(46.56, 7.59)}  & (46.56, 7.59)  & \multicolumn{1}{l}{(72.52, 4.42)} & \multicolumn{1}{l}{(69.94, 5.35)} & \multicolumn{1}{l}{(60.99, 8.01)}  & (53.12, 9.94)  \\ 
		& 0.2                          & \multicolumn{1}{l}{(29.38, 3.42)}  & \multicolumn{1}{l}{(29.38, 3.42)}  & \multicolumn{1}{l}{(29.38, 3.42)}  & (29.38, 3.42)  & \multicolumn{1}{l}{(66.79, 3.34)} & \multicolumn{1}{l}{(64.67, 4.05)} & \multicolumn{1}{l}{(56.89, 6.07)}  & (49.79, 7.55)  \\ 
		& 0.4                          & \multicolumn{1}{l}{(9.05, 0.71)}   & \multicolumn{1}{l}{(9.05, 0.71)}   & \multicolumn{1}{l}{(9.05, 0.71)}   & (9.05, 0.71)   & \multicolumn{1}{l}{(52.03, 2.02)} & \multicolumn{1}{l}{(50.93, 2.44)} & \multicolumn{1}{l}{(45.83, 3.64)}  & (40.61, 4.51)  \\ 
		& 0.6                          & \multicolumn{1}{l}{(2.32, 0.15)}   & \multicolumn{1}{l}{(2.32, 0.15)}   & \multicolumn{1}{l}{(2.32, 0.15)}   & (2.32, 0.15)   & \multicolumn{1}{l}{(37.76, 1.23)} & \multicolumn{1}{l}{(37.46, 1.48)} & \multicolumn{1}{l}{(34.68, 2.21)}  & (31.13, 2.73)  \\ 
		& 0.9                          & \multicolumn{1}{l}{(0.29, 0.01)}   & \multicolumn{1}{l}{(0.29, 0.01)}   & \multicolumn{1}{l}{(0.29, 0.01)}   & (0.29, 0.01)   & \multicolumn{1}{l}{(21.52, 0.56)} & \multicolumn{1}{l}{(21.76, 0.67)} & \multicolumn{1}{l}{(21.04, 1.01)}  & (19.36, 1.25)  \\ 
		& 1.2                          & \multicolumn{1}{l}{(0.02, 0.00)}   & \multicolumn{1}{l}{(0.02, 0.00)}   & \multicolumn{1}{l}{(0.02, 0.00)}   & (0.02, 0.00)   & \multicolumn{1}{l}{(11.59, 0.25)} & \multicolumn{1}{l}{(11.97, 0.31)} & \multicolumn{1}{l}{(12.13, 0.46)}  & (11.45, 0.57)  \\ 
		& 1.4                          & \multicolumn{1}{l}{(0.00, 0.00)}   & \multicolumn{1}{l}{(0.00, 0.00)}   & \multicolumn{1}{l}{(0.00, 0.00)}   & (0.00, 0.00)   & \multicolumn{1}{l}{(7.50, 0.15)}  & \multicolumn{1}{l}{(7.84, 0.18)}  & \multicolumn{1}{l}{(8.18, 0.27)}   & (7.85, 0.34)   \\ \hline\hline
		\multirow{7}{*}{$(15,12)$}  & 0.1                          & \multicolumn{1}{l}{(26.55, 6.37)}  & \multicolumn{1}{l}{(26.96, 6.92)}  & \multicolumn{1}{l}{(27.81, 8.30)}  & (28.24, 9.18)  & \multicolumn{1}{l}{(60.51, 6.57)} & \multicolumn{1}{l}{(59.16, 7.50)} & \multicolumn{1}{l}{(55.38, 9.97)}  & (52.75, 11.63) \\ 
		& 0.2                          & \multicolumn{1}{l}{(8.96, 1.45)}   & \multicolumn{1}{l}{(9.27, 1.57)}   & \multicolumn{1}{l}{(9.96, 1.87)}   & (10.33, 2.06)  & \multicolumn{1}{l}{(48.53, 3.93)} & \multicolumn{1}{l}{(47.69, 4.49)} & \multicolumn{1}{l}{(45.12, 5.99)}  & (43.20, 6.99)  \\ 
		& 0.4                          & \multicolumn{1}{l}{(0.66, 0.08)}   & \multicolumn{1}{l}{(0.71, 0.09)}   & \multicolumn{1}{l}{(0.82, 0.10)}   & (0.89, 0.11)   & \multicolumn{1}{l}{(26.20, 1.40)} & \multicolumn{1}{l}{(26.07, 1.60)} & \multicolumn{1}{l}{(25.30, 2.14)}  & (24.54, 2.50)  \\ 
		& 0.6                          & \multicolumn{1}{l}{(0.00, 0.01)}   & \multicolumn{1}{l}{(0.00, 0.01)}   & \multicolumn{1}{l}{(0.00, 0.01)}   & (0.00, 0.01)   & \multicolumn{1}{l}{(12.41, 0.54)} & \multicolumn{1}{l}{(12.50, 0.62)} & \multicolumn{1}{l}{(12.43, 0.82)}  & (12.20, 0.95)  \\ 
		& 0.9                          & \multicolumn{1}{l}{(0.00, 0.00)}   & \multicolumn{1}{l}{(0.00, 0.00)}   & \multicolumn{1}{l}{(0.00, 0.00)}   & (0.00, 0.00)   & \multicolumn{1}{l}{(3.25, 0.14)}  & \multicolumn{1}{l}{(3.33, 0.15)}  & \multicolumn{1}{l}{(3.44, 0.20)}   & (3.45, 0.23)   \\ 
		& 1.2                          & \multicolumn{1}{l}{(0.00, 0.00)}   & \multicolumn{1}{l}{(0.00, 0.00)}   & \multicolumn{1}{l}{(0.00, 0.00)}   & (0.00, 0.00)   & \multicolumn{1}{l}{(0.79, 0.04)}  & \multicolumn{1}{l}{(0.83, 0.04)}  & \multicolumn{1}{l}{(0.88, 0.05)}   & (0.90, 0.06)   \\ 
		& 1.4                          & \multicolumn{1}{l}{(0.00, 0.00)}   & \multicolumn{1}{l}{(0.00, 0.00)}   & \multicolumn{1}{l}{(0.00, 0.00)}   & (0.00, 0.00)   & \multicolumn{1}{l}{(0.23, 0.02)}  & \multicolumn{1}{l}{(0.25, 0.02)}  & \multicolumn{1}{l}{(0.26, 0.03)}   & (0.27, 0.03)   \\ \hline \hline
	\end{tabular}}
\end{table}

\begin{table}[H]
		\caption{PRIs of the estimators $\delta_{\psi_{01}}$ and $\delta_{\psi_{02}}$ with respect to the BAEE when $\sigma_i$'s are known }
		\vspace{0.2cm}
		\label{tab7}
	\centering
	\resizebox{18cm}{10cm}{

			\begin{tabular}{llllllllll}
				\hline \hline
				\multirow{3}{*}{$(n_1,n_2)$} & \multirow{3}{*}{$\mu_2-\mu_1$} & \multicolumn{4}{c}{$(\sigma_1,\sigma_2)=(1,1.5)$}                                                                             & \multicolumn{4}{c}{$(\sigma_1,\sigma_2)=(3,2)$}                                                                               \\ \cline{3-10} 
				&                              & \multicolumn{1}{l}{$p=-2.5$}      & \multicolumn{1}{l}{$p=-2$}        & \multicolumn{1}{l}{$p=2$}          & $p=2.5$        & \multicolumn{1}{l}{$p=-2.5$}      & \multicolumn{1}{l}{$p=-2$}        & \multicolumn{1}{l}{$p=2$}          & $p=2.5$        \\ \cline{3-10} 
				&                              & \multicolumn{1}{l}{$(\delta_{\psi_{01}},\delta_{\psi_{02}})$}             & \multicolumn{1}{l}{$(\delta_{\psi_{01}},\delta_{\psi_{02}})$}             & \multicolumn{1}{l}{$(\delta_{\psi_{01}},\delta_{\psi_{02}})$}             & {$(\delta_{\psi_{01}},\delta_{\psi_{02}})$}              & \multicolumn{1}{l}{$(\delta_{\psi_{01}},\delta_{\psi_{02}})$}            & \multicolumn{1}{l}{$(\delta_{\psi_{01}},\delta_{\psi_{02}})$}             & \multicolumn{1}{l}{$(\delta_{\psi_{01}},\delta_{\psi_{02}})$}             & {$(\delta_{\psi_{01}},\delta_{\psi_{02}})$}             \\ \hline \hline
				\multirow{7}{*}{$(8,8)$}    & 0.1                          & \multicolumn{1}{l}{(40.86, 6.74)} & \multicolumn{1}{l}{(41.35, 7.42)} & \multicolumn{1}{l}{(45.56, 21.12)} & (45.56, 21.12) & \multicolumn{1}{l}{(73.65, 2.93)} & \multicolumn{1}{l}{(72.76, 3.55)} & \multicolumn{1}{l}{(23.08, 24.32)} & (3.67, 33.59)  \\ 
				& 0.2                          & \multicolumn{1}{l}{(24.62, 3.03)} & \multicolumn{1}{l}{(25.32, 3.33)} & \multicolumn{1}{l}{(32.57, 9.88)}  & (32.57, 9.88)  & \multicolumn{1}{l}{(67.02, 2.23)} & \multicolumn{1}{l}{(66.47, 2.70)} & \multicolumn{1}{l}{(21.80, 19.02)} & (3.48, 27.14)  \\ 
				& 0.4                          & \multicolumn{1}{l}{(6.84, 0.66)}  & \multicolumn{1}{l}{(7.25, 0.72)}  & \multicolumn{1}{l}{(13.33, 2.05)}  & (13.33, 2.05)  & \multicolumn{1}{l}{(50.46, 1.33)} & \multicolumn{1}{l}{(50.54, 1.61)} & \multicolumn{1}{l}{(18.21, 11.14)} & (2.93, 15.88)  \\ 
				& 0.6                          & \multicolumn{1}{l}{(1.56, 0.14)}  & \multicolumn{1}{l}{(1.70, 0.15)}  & \multicolumn{1}{l}{(4.58, 0.44)}   & (4.58, 0.44)   & \multicolumn{1}{l}{(35.35, 0.75)} & \multicolumn{1}{l}{(35.81, 0.91)} & \multicolumn{1}{l}{(14.30, 6.53)}  & (2.34, 9.34)   \\ 
				& 0.9                          & \multicolumn{1}{l}{(0.15, 0.01)}  & \multicolumn{1}{l}{(0.17, 0.01)}  & \multicolumn{1}{l}{(0.87, 0.03)}   & (0.87, 0.03)   & \multicolumn{1}{l}{(19.27, 0.33)} & \multicolumn{1}{l}{(19.82, 0.40)} & \multicolumn{1}{l}{(9.29, 2.88)}   & (1.56, 4.12)   \\ 
				& 1.2                          & \multicolumn{1}{l}{(0.02, 0.00)}  & \multicolumn{1}{l}{(0.02, 0.00)}  & \multicolumn{1}{l}{(0.14, 0.00)}   & (0.14, 0.00)   & \multicolumn{1}{l}{(9.77, 0.19)}  & \multicolumn{1}{l}{(10.20, 0.22)} & \multicolumn{1}{l}{(5.73, 1.39)}   & (1.00, 1.96)   \\ 
				& 1.4                          & \multicolumn{1}{l}{(0.00, 0.00)}  & \multicolumn{1}{l}{(0.00, 0.00)}  & \multicolumn{1}{l}{(0.00, 0.00)}   & (0.00, 0.00)   & \multicolumn{1}{l}{(6.07, 0.11)}  & \multicolumn{1}{l}{(6.40, 0.13)}  & \multicolumn{1}{l}{(4.11, 0.86)}   & (0.71, 1.20)   \\ \hline \hline
				\multirow{7}{*}{$(9,10)$}   & 0.1                          & \multicolumn{1}{l}{(40.52, 6.10)} & \multicolumn{1}{l}{(40.99, 6.63)} & \multicolumn{1}{l}{(44.79, 15.77)} & (45.18, 18.27) & \multicolumn{1}{l}{(74.06, 3.15)} & \multicolumn{1}{l}{(73.09, 3.71)} & \multicolumn{1}{l}{(30.67, 18.06)} & (13.53, 23.24) \\ 
				& 0.2                          & \multicolumn{1}{l}{(22.14, 2.44)} & \multicolumn{1}{l}{(22.77, 2.65)} & \multicolumn{1}{l}{(29.08, 6.32)}  & (30.00, 7.37)  & \multicolumn{1}{l}{(66.13, 2.25)} & \multicolumn{1}{l}{(65.52, 2.66)} & \multicolumn{1}{l}{(28.47, 13.28)} & (12.60, 17.23) \\ 
				& 0.4                          & \multicolumn{1}{l}{(4.73, 0.43)}  & \multicolumn{1}{l}{(5.01, 0.46)}  & \multicolumn{1}{l}{(8.98, 1.07)}   & (9.77, 1.24)   & \multicolumn{1}{l}{(47.19, 1.23)} & \multicolumn{1}{l}{(47.25, 1.46)} & \multicolumn{1}{l}{(22.56, 7.27)}  & (10.04, 9.42)  \\ 
				& 0.6                          & \multicolumn{1}{l}{(0.78, 0.07)}  & \multicolumn{1}{l}{(0.85, 0.07)}  & \multicolumn{1}{l}{(2.21, 0.18)}   & (2.55, 0.20)   & \multicolumn{1}{l}{(31.12, 0.69)} & \multicolumn{1}{l}{(31.50, 0.82)} & \multicolumn{1}{l}{(16.58, 4.05)}  & (7.43, 5.24)   \\ 
				& 0.9                          & \multicolumn{1}{l}{(0.06, 0.01)}  & \multicolumn{1}{l}{(0.06, 0.01)}  & \multicolumn{1}{l}{(0.24, 0.01)}   & (0.30, 0.01)   & \multicolumn{1}{l}{(15.16, 0.25)} & \multicolumn{1}{l}{(15.58, 0.30)} & \multicolumn{1}{l}{(9.53, 1.60)}   & (4.32, 2.09)   \\ 
				& 1.2                          & \multicolumn{1}{l}{(0.00, 0.00)}  & \multicolumn{1}{l}{(0.00, 0.00)}  & \multicolumn{1}{l}{(0.00, 0.00)}   & (0.00, 0.00)   & \multicolumn{1}{l}{(6.78, 0.13)}  & \multicolumn{1}{l}{(7.07, 0.15)}  & \multicolumn{1}{l}{(4.98, 0.70)}   & (2.30, 0.91)   \\ 
				& 1.4                          & \multicolumn{1}{l}{(0.00, 0.00)}  & \multicolumn{1}{l}{(0.00, 0.00)}  & \multicolumn{1}{l}{(0.00, 0.00)}   & (0.00, 0.00)   & \multicolumn{1}{l}{(3.79, 0.06)}  & \multicolumn{1}{l}{(3.99, 0.08)}  & \multicolumn{1}{l}{(3.14, 0.39)}   & (1.49, 0.50)   \\ \hline \hline
				\multirow{7}{*}{$(12,8)$}   & 0.1                          & \multicolumn{1}{l}{(28.70, 7.45)} & \multicolumn{1}{l}{(29.14, 8.23)} & \multicolumn{1}{l}{(33.11, 23.24)} & (33.64, 27.52) & \multicolumn{1}{l}{(63.86, 4.99)} & \multicolumn{1}{l}{(63.17, 5.81)} & \multicolumn{1}{l}{(40.02, 29.73)} & (31.51, 39.40) \\ 
				& 0.2                          & \multicolumn{1}{l}{(12.53, 2.32)} & \multicolumn{1}{l}{(12.94, 2.56)} & \multicolumn{1}{l}{(17.24, 8.34)}  & (17.91, 10.72) & \multicolumn{1}{l}{(53.96, 3.31)} & \multicolumn{1}{l}{(53.63, 3.86)} & \multicolumn{1}{l}{(35.29, 21.64)} & (27.93, 30.27) \\ 
				& 0.4                          & \multicolumn{1}{l}{(1.73, 0.26)}  & \multicolumn{1}{l}{(1.83, 0.28)}  & \multicolumn{1}{l}{(3.25, 0.83)}   & (3.52, 1.05)   & \multicolumn{1}{l}{(33.29, 1.55)} & \multicolumn{1}{l}{(33.46, 1.80)} & \multicolumn{1}{l}{(24.38, 9.76)}  & (19.55, 13.77) \\ 
				& 0.6                          & \multicolumn{1}{l}{(0.08, 0.03)}  & \multicolumn{1}{l}{(0.09, 0.04)}  & \multicolumn{1}{l}{(0.20, 0.11)}   & (0.22, 0.13)   & \multicolumn{1}{l}{(18.57, 0.73)} & \multicolumn{1}{l}{(18.89, 0.84)} & \multicolumn{1}{l}{(15.46, 4.56)}  & (12.60, 6.40)  \\ 
				& 0.9                          & \multicolumn{1}{l}{(0.00, 0.00)}  & \multicolumn{1}{l}{(0.00, 0.00)}  & \multicolumn{1}{l}{(0.00, 0.00)}   & (0.00, 0.00)   & \multicolumn{1}{l}{(7.10, 0.24)}  & \multicolumn{1}{l}{(7.33, 0.28)}  & \multicolumn{1}{l}{(6.90, 1.41)}   & (5.72, 1.98)   \\ 
				& 1.2                          & \multicolumn{1}{l}{(0.00, 0.00)}  & \multicolumn{1}{l}{(0.00, 0.00)}  & \multicolumn{1}{l}{(0.00, 0.00)}   & (0.00, 0.00)   & \multicolumn{1}{l}{(2.46, 0.09)}  & \multicolumn{1}{l}{(2.57, 0.10)}  & \multicolumn{1}{l}{(2.63, 0.49)}   & (2.14, 0.67)   \\ 
				& 1.4                          & \multicolumn{1}{l}{(0.00, 0.00)}  & \multicolumn{1}{l}{(0.00, 0.00)}  & \multicolumn{1}{l}{(0.00, 0.00)}   & (0.00, 0.00)   & \multicolumn{1}{l}{(1.05, 0.05)}  & \multicolumn{1}{l}{(1.10, 0.06)}  & \multicolumn{1}{l}{(1.10, 0.24)}   & (0.87, 0.32)   \\ \hline \hline
				\multirow{7}{*}{$(14,15)$}  & 0.1                          & \multicolumn{1}{l}{(29.83, 4.48)} & \multicolumn{1}{l}{(30.25, 4.75)} & \multicolumn{1}{l}{(34.04, 8.35)}  & (34.58, 9.11)  & \multicolumn{1}{l}{(68.04, 3.84)} & \multicolumn{1}{l}{(67.19, 4.29)} & \multicolumn{1}{l}{(46.43, 11.77)} & (39.93, 13.66) \\ 
				& 0.2                          & \multicolumn{1}{l}{(10.44, 1.11)} & \multicolumn{1}{l}{(10.79, 1.17)} & \multicolumn{1}{l}{(14.63, 2.05)}  & (15.29, 2.23)  & \multicolumn{1}{l}{(54.78, 2.39)} & \multicolumn{1}{l}{(54.36, 2.67)} & \multicolumn{1}{l}{(39.21, 7.34)}  & (33.91, 8.52)  \\ 
				& 0.4                          & \multicolumn{1}{l}{(1.06, 0.07)}  & \multicolumn{1}{l}{(1.13, 0.08)}  & \multicolumn{1}{l}{(2.26, 0.13)}   & (2.52, 0.14)   & \multicolumn{1}{l}{(29.66, 0.94)} & \multicolumn{1}{l}{(29.79, 1.05)} & \multicolumn{1}{l}{(24.08, 2.91)}  & (21.19, 3.38)  \\ 
				& 0.6                          & \multicolumn{1}{l}{(0.10, 0.00)}  & \multicolumn{1}{l}{(0.11, 0.00)}  & \multicolumn{1}{l}{(0.31, 0.01)}   & (0.37, 0.01)   & \multicolumn{1}{l}{(14.34, 0.34)} & \multicolumn{1}{l}{(14.59, 0.38)} & \multicolumn{1}{l}{(13.45, 1.06)}  & (12.09, 1.24)  \\ 
				& 0.9                          & \multicolumn{1}{l}{(0.00, 0.00)}  & \multicolumn{1}{l}{(0.00, 0.00)}  & \multicolumn{1}{l}{(0.00, 0.00)}   & (0.00, 0.00)   & \multicolumn{1}{l}{(4.51, 0.08)}  & \multicolumn{1}{l}{(4.68, 0.09)}  & \multicolumn{1}{l}{(5.33, 0.25)}   & (4.96, 0.29)   \\ 
				& 1.2                          & \multicolumn{1}{l}{(0.00, 0.00)}  & \multicolumn{1}{l}{(0.00, 0.00)}  & \multicolumn{1}{l}{(0.00, 0.00)}   & (0.00, 0.00)   & \multicolumn{1}{l}{(1.40, 0.02)}  & \multicolumn{1}{l}{(1.48, 0.03)}  & \multicolumn{1}{l}{(2.04, 0.06)}   & (1.94, 0.07)   \\ 
				& 1.4                          & \multicolumn{1}{l}{(0.00, 0.00)}  & \multicolumn{1}{l}{(0.00, 0.00)}  & \multicolumn{1}{l}{(0.00, 0.00)}   & (0.00, 0.00)   & \multicolumn{1}{l}{(0.59, 0.01)}  & \multicolumn{1}{l}{(0.63, 0.02)}  & \multicolumn{1}{l}{(1.00, 0.04)}   & (0.99, 0.04)   \\ \hline \hline
				\multirow{7}{*}{$(16,13)$}  & 0.1                          & \multicolumn{1}{l}{(23.69, 4.91)} & \multicolumn{1}{l}{(24.05, 5.25)} & \multicolumn{1}{l}{(27.47, 10.29)} & (27.96, 11.49) & \multicolumn{1}{l}{(62.76, 4.85)} & \multicolumn{1}{l}{(62.07, 5.41)} & \multicolumn{1}{l}{(46.31, 15.52)} & (41.70, 18.34) \\ 
				& 0.2                          & \multicolumn{1}{l}{(6.81, 1.01)}  & \multicolumn{1}{l}{(7.05, 1.07)}  & \multicolumn{1}{l}{(9.71, 2.08)}   & (10.16, 2.33)  & \multicolumn{1}{l}{(48.27, 2.85)} & \multicolumn{1}{l}{(47.98, 3.18)} & \multicolumn{1}{l}{(37.42, 9.19)}  & (33.91, 10.85) \\ 
				& 0.4                          & \multicolumn{1}{l}{(0.54, 0.05)}  & \multicolumn{1}{l}{(0.58, 0.05)}  & \multicolumn{1}{l}{(1.19, 0.09)}   & (1.33, 0.10)   & \multicolumn{1}{l}{(23.27, 0.99)} & \multicolumn{1}{l}{(23.41, 1.10)} & \multicolumn{1}{l}{(20.41, 3.18)}  & (18.80, 3.76)  \\ 
				& 0.6                          & \multicolumn{1}{l}{(0.01, 0.00)}  & \multicolumn{1}{l}{(0.01, 0.00)}  & \multicolumn{1}{l}{(0.04, 0.00)}   & (0.05, 0.00)   & \multicolumn{1}{l}{(9.85, 0.33)}  & \multicolumn{1}{l}{(10.05, 0.37)} & \multicolumn{1}{l}{(10.11, 1.06)}  & (9.52, 1.25)   \\ 
				& 0.9                          & \multicolumn{1}{l}{(0.00, 0.00)}  & \multicolumn{1}{l}{(0.00, 0.00)}  & \multicolumn{1}{l}{(0.00, 0.00)}   & (0.00, 0.00)   & \multicolumn{1}{l}{(2.53, 0.06)}  & \multicolumn{1}{l}{(2.63, 0.07)}  & \multicolumn{1}{l}{(3.27, 0.20)}   & (3.20, 0.24)   \\ 
				& 1.2                          & \multicolumn{1}{l}{(0.00, 0.00)}  & \multicolumn{1}{l}{(0.00, 0.00)}  & \multicolumn{1}{l}{(0.00, 0.00)}   & (0.00, 0.00)   & \multicolumn{1}{l}{(0.70, 0.01)}  & \multicolumn{1}{l}{(0.74, 0.02)}  & \multicolumn{1}{l}{(1.16, 0.04)}   & (1.18, 0.04)   \\ 
				& 1.4                          & \multicolumn{1}{l}{(0.00, 0.00)}  & \multicolumn{1}{l}{(0.00, 0.00)}  & \multicolumn{1}{l}{(0.00, 0.00)}   & (0.00, 0.00)   & \multicolumn{1}{l}{(0.28, 0.01)}  & \multicolumn{1}{l}{(0.30, 0.01)}  & \multicolumn{1}{l}{(0.52, 0.02)}   & (0.54, 0.02)   \\ \hline \hline
			\end{tabular}}
	\end{table}

   \begin{table}[H]
   		\caption{PRIs of the estimator $\delta_{1RMLI}$ with respect to RMLE when $\sigma_i$'s are known. }
   		\vspace{0.2cm}
   \label{tab1RML}
   	\centering
   	\resizebox{16cm}{10.5cm}{
   	
   		\begin{tabular}{llllllllll}
   			\hline
   			\multirow{3}{*}{$(n_1,n_2)$} & \multirow{3}{*}{$\mu_2-\mu_1$} & \multicolumn{4}{c}{$(\sigma_1,\sigma_2)=(1,1.5)$}                                                 & \multicolumn{4}{c}{$(\sigma_1,\sigma_2)=(3,2)$}                                                   \\ \cline{3-10} 
   			&                                & \multicolumn{1}{l}{$p=-1$} & \multicolumn{1}{l}{$p=-0.5$} & \multicolumn{1}{l}{$p=0.5$} & $p=1$ & \multicolumn{1}{l}{$p=-1$} & \multicolumn{1}{l}{$p=-0.5$} & \multicolumn{1}{l}{$p=0.5$} & $p=1$ \\ \cline{3-10} 
   			&                                & \multicolumn{1}{l}{$\delta_{1RMLI}$}     & \multicolumn{1}{l}{$\delta_{1RMLI}$}       & \multicolumn{1}{l}{$\delta_{1RMLI}$}      & $\delta_{1RMLI}$    & \multicolumn{1}{l}{$\delta_{1RMLI}$}     & \multicolumn{1}{l}{$\delta_{1RMLI}$}       & \multicolumn{1}{l}{$\delta_{1RMLI}$}      & $\delta_{1RMLI}$    \\ \hline \hline
   			\multirow{7}{*}{$(5,5)$}     & 0.1                            & \multicolumn{1}{l}{53.73}  & \multicolumn{1}{l}{54.77}    & \multicolumn{1}{l}{57.03}   & 58.25 & \multicolumn{1}{l}{55.72}  & \multicolumn{1}{l}{56.36}    & \multicolumn{1}{l}{57.71}   & 58.43 \\ 
   			& 0.2                            & \multicolumn{1}{l}{51.08}  & \multicolumn{1}{l}{52.31}    & \multicolumn{1}{l}{55.00}   & 56.46 & \multicolumn{1}{l}{56.34}  & \multicolumn{1}{l}{57.09}    & \multicolumn{1}{l}{58.66}   & 59.49 \\ 
   			& 0.4                            & \multicolumn{1}{l}{45.35}  & \multicolumn{1}{l}{46.65}    & \multicolumn{1}{l}{49.54}   & 51.16 & \multicolumn{1}{l}{52.05}  & \multicolumn{1}{l}{52.98}    & \multicolumn{1}{l}{54.94}   & 55.98 \\ 
   			& 0.6                            & \multicolumn{1}{l}{42.31}  & \multicolumn{1}{l}{43.52}    & \multicolumn{1}{l}{46.27}   & 47.84 & \multicolumn{1}{l}{47.08}  & \multicolumn{1}{l}{48.09}    & \multicolumn{1}{l}{50.23}   & 51.37 \\ 
   			& 0.9                            & \multicolumn{1}{l}{40.56}  & \multicolumn{1}{l}{41.63}    & \multicolumn{1}{l}{44.09}   & 45.51 & \multicolumn{1}{l}{41.43}  & \multicolumn{1}{l}{42.44}    & \multicolumn{1}{l}{44.63}   & 45.81 \\ 
   			& 1.2                            & \multicolumn{1}{l}{40.04}  & \multicolumn{1}{l}{41.04}    & \multicolumn{1}{l}{43.31}   & 44.62 & \multicolumn{1}{l}{37.69}  & \multicolumn{1}{l}{38.64}    & \multicolumn{1}{l}{40.74}   & 41.91 \\ 
   			& 1.4                            & \multicolumn{1}{l}{39.94}  & \multicolumn{1}{l}{40.92}    & \multicolumn{1}{l}{43.14}   & 44.41 & \multicolumn{1}{l}{35.97}  & \multicolumn{1}{l}{36.87}    & \multicolumn{1}{l}{38.90}   & 40.03 \\ \hline \hline
   			\multirow{7}{*}{$(5,7)$}     & 0.1                            & \multicolumn{1}{l}{54.61}  & \multicolumn{1}{l}{55.57}    & \multicolumn{1}{l}{57.61}   & 58.70 & \multicolumn{1}{l}{57.33}  & \multicolumn{1}{l}{57.87}    & \multicolumn{1}{l}{58.98}   & 58.98 \\ 
   			& 0.2                            & \multicolumn{1}{l}{50.34}  & \multicolumn{1}{l}{51.49}    & \multicolumn{1}{l}{53.98}   & 55.32 & \multicolumn{1}{l}{56.56}  & \multicolumn{1}{l}{57.22}    & \multicolumn{1}{l}{58.60}   & 58.60 \\ 
   			& 0.4                            & \multicolumn{1}{l}{42.97}  & \multicolumn{1}{l}{44.15}    & \multicolumn{1}{l}{46.77}   & 48.22 & \multicolumn{1}{l}{49.59}  & \multicolumn{1}{l}{50.42}    & \multicolumn{1}{l}{52.15}   & 52.15 \\ 
   			& 0.6                            & \multicolumn{1}{l}{39.34}  & \multicolumn{1}{l}{40.40}    & \multicolumn{1}{l}{42.80}   & 44.17 & \multicolumn{1}{l}{43.30}  & \multicolumn{1}{l}{44.17}    & \multicolumn{1}{l}{46.01}   & 46.01 \\ 
   			& 0.9                            & \multicolumn{1}{l}{37.30}  & \multicolumn{1}{l}{38.20}    & \multicolumn{1}{l}{40.27}   & 41.46 & \multicolumn{1}{l}{36.92}  & \multicolumn{1}{l}{37.76}    & \multicolumn{1}{l}{39.57}   & 39.57 \\ 
   			& 1.2                            & \multicolumn{1}{l}{36.73}  & \multicolumn{1}{l}{37.56}    & \multicolumn{1}{l}{39.42}   & 40.48 & \multicolumn{1}{l}{33.03}  & \multicolumn{1}{l}{33.81}    & \multicolumn{1}{l}{35.52}   & 35.52 \\ 
   			& 1.4                            & \multicolumn{1}{l}{36.60}  & \multicolumn{1}{l}{37.40}    & \multicolumn{1}{l}{39.19}   & 40.20 & \multicolumn{1}{l}{31.29}  & \multicolumn{1}{l}{32.02}    & \multicolumn{1}{l}{33.65}   & 33.65 \\ \hline \hline
   			\multirow{7}{*}{$(7,6)$}     & 0.1                            & \multicolumn{1}{l}{53.32}  & \multicolumn{1}{l}{53.32}    & \multicolumn{1}{l}{55.94}   & 56.88 & \multicolumn{1}{l}{56.37}  & \multicolumn{1}{l}{56.90}    & \multicolumn{1}{l}{57.99}   & 58.57 \\ 
   			& 0.2                            & \multicolumn{1}{l}{49.67}  & \multicolumn{1}{l}{49.67}    & \multicolumn{1}{l}{52.70}   & 53.81 & \multicolumn{1}{l}{55.73}  & \multicolumn{1}{l}{56.37}    & \multicolumn{1}{l}{57.68}   & 58.37 \\ 
   			& 0.4                            & \multicolumn{1}{l}{44.46}  & \multicolumn{1}{l}{44.46}    & \multicolumn{1}{l}{47.45}   & 48.58 & \multicolumn{1}{l}{55.73}  & \multicolumn{1}{l}{50.63}    & \multicolumn{1}{l}{52.23}   & 53.07 \\
   			& 0.6                            & \multicolumn{1}{l}{42.57}  & \multicolumn{1}{l}{45.31}   & 46.36 & \multicolumn{1}{l}{48.04}  &45.08& \multicolumn{1}{l}{44.70}    & \multicolumn{1}{l}{47.16}   & 48.04 \\ 
   			& 0.9                            & \multicolumn{1}{l}{41.89}  & \multicolumn{1}{l}{41.89}    & \multicolumn{1}{l}{44.42}   & 45.38 & \multicolumn{1}{l}{42.90}  & \multicolumn{1}{l}{40.42}    & \multicolumn{1}{l}{42.03}   & 42.90 \\ 
   			& 1.2                            & \multicolumn{1}{l}{41.77}  & \multicolumn{1}{l}{41.77}    & \multicolumn{1}{l}{44.23}   & 45.16 & \multicolumn{1}{l}{39.75}  & \multicolumn{1}{l}{37.45}    & \multicolumn{1}{l}{38.94}   & 39.75 \\ 
   			& 1.4                            & \multicolumn{1}{l}{41.75}  & \multicolumn{1}{l}{41.75}    & \multicolumn{1}{l}{44.19}   & 45.11 & \multicolumn{1}{l}{38.40}  & \multicolumn{1}{l}{36.24}    & \multicolumn{1}{l}{37.63}   & 38.40 \\ \hline \hline
   			\multirow{7}{*}{$(8,10$}     & 0.1                            & \multicolumn{1}{l}{53.32}  & \multicolumn{1}{l}{54.05}    & \multicolumn{1}{l}{55.56}   & 56.35 & \multicolumn{1}{l}{57.61}  & \multicolumn{1}{l}{57.61}    & \multicolumn{1}{l}{57.61}   & 57.61 \\ 
   			& 0.2                            & \multicolumn{1}{l}{47.08}  & \multicolumn{1}{l}{47.89}    & \multicolumn{1}{l}{49.60}   & 50.50 & \multicolumn{1}{l}{53.80}  & \multicolumn{1}{l}{53.80}    & \multicolumn{1}{l}{53.80}   & 53.80 \\ 
   			& 0.4                            & \multicolumn{1}{l}{40.78}  & \multicolumn{1}{l}{41.50}    & \multicolumn{1}{l}{43.05}   & 43.89 & \multicolumn{1}{l}{44.39}  & \multicolumn{1}{l}{44.39}    & \multicolumn{1}{l}{44.39}   & 44.39 \\ 
   			& 0.6                            & \multicolumn{1}{l}{38.93}  & \multicolumn{1}{l}{39.55}    & \multicolumn{1}{l}{40.90}   & 41.64 & \multicolumn{1}{l}{38.28}  & \multicolumn{1}{l}{38.28}    & \multicolumn{1}{l}{38.28}   & 38.28 \\ 
   			& 0.9                            & \multicolumn{1}{l}{38.36}  & \multicolumn{1}{l}{38.93}    & \multicolumn{1}{l}{40.14}   & 40.79 & \multicolumn{1}{l}{33.20}  & \multicolumn{1}{l}{33.20}    & \multicolumn{1}{l}{33.20}   & 33.20 \\ 
   			& 1.2                            & \multicolumn{1}{l}{38.29}  & \multicolumn{1}{l}{38.84}    & \multicolumn{1}{l}{40.02}   & 40.64 & \multicolumn{1}{l}{30.63}  & \multicolumn{1}{l}{30.63}    & \multicolumn{1}{l}{30.63}   & 30.63 \\ 
   			& 1.4                            & \multicolumn{1}{l}{38.29}  & \multicolumn{1}{l}{38.83}    & \multicolumn{1}{l}{40.00}   & 40.62 & \multicolumn{1}{l}{29.67}  & \multicolumn{1}{l}{29.67}    & \multicolumn{1}{l}{29.67}   & 29.67 \\ \hline \hline
   			\multirow{7}{*}{$(15,12$}    & 0.1                            & \multicolumn{1}{l}{50.31}  & \multicolumn{1}{l}{50.79}    & \multicolumn{1}{l}{51.77}   & 51.77 & \multicolumn{1}{l}{55.99}  & \multicolumn{1}{l}{56.30}    & \multicolumn{1}{l}{56.30}   & 57.24 \\ 
   			& 0.2                            & \multicolumn{1}{l}{45.45}  & \multicolumn{1}{l}{45.92}    & \multicolumn{1}{l}{46.88}   & 46.88 & \multicolumn{1}{l}{50.32}  & \multicolumn{1}{l}{50.69}    & \multicolumn{1}{l}{50.69}   & 51.84 \\ 
   			& 0.4                            & \multicolumn{1}{l}{43.20}  & \multicolumn{1}{l}{43.60}    & \multicolumn{1}{l}{44.44}   & 44.44 & \multicolumn{1}{l}{41.78}  & \multicolumn{1}{l}{42.16}    & \multicolumn{1}{l}{42.16}   & 43.34 \\
   			& 0.6                            & \multicolumn{1}{l}{43.01}  & \multicolumn{1}{l}{43.40}    & \multicolumn{1}{l}{44.21}   & 44.21 & \multicolumn{1}{l}{37.76}  & \multicolumn{1}{l}{38.09}    & \multicolumn{1}{l}{38.09}   & 39.15 \\ 
   			& 0.9                            & \multicolumn{1}{l}{43.01}  & \multicolumn{1}{l}{43.40}    & \multicolumn{1}{l}{44.21}   & 44.21 & \multicolumn{1}{l}{35.39}  & \multicolumn{1}{l}{35.67}    & \multicolumn{1}{l}{35.67}   & 36.56 \\ 
   			& 1.2                            & \multicolumn{1}{l}{43.01}  & \multicolumn{1}{l}{43.40}    & \multicolumn{1}{l}{44.21}   & 44.21 & \multicolumn{1}{l}{34.74}  & \multicolumn{1}{l}{34.99}    & \multicolumn{1}{l}{34.99}   & 35.78 \\ 
   			& 1.4                            & \multicolumn{1}{l}{43.01}  & \multicolumn{1}{l}{43.40}    & \multicolumn{1}{l}{44.21}   & 44.21 & \multicolumn{1}{l}{34.60}  & \multicolumn{1}{l}{34.84}    & \multicolumn{1}{l}{34.84}   & 35.60 \\ \hline \hline
   		\end{tabular}}
   	\end{table}
   	
   	\begin{table}[H]
   			\caption{PRIs of the estimator $\delta_{1RMLI}$ with respect to RMLE when $\sigma_i$'s are known.}
   			\vspace{0.2cm}
   		 \label{tab2RML}
   		\centering
   		\resizebox{16cm}{10.5cm}{
   			
   		\begin{tabular}{llllllllll}
   			\hline \hline
   			\multirow{3}{*}{$(n_1,n_2)$} & \multirow{3}{*}{$\mu_2-\mu_1$} & \multicolumn{4}{c}{$(\sigma_1,\sigma_2)=(1,1.5)$}                                                  & \multicolumn{4}{c}{$(\sigma_1,\sigma_2)=(3,2))$}                                                    \\ \cline{3-10} 
   			&                                & \multicolumn{1}{l}{$p=-2.5$} & \multicolumn{1}{l}{$p=-2$} & \multicolumn{1}{l}{$p=2$} & $p=2.5$ & \multicolumn{1}{l}{$p=-2.5$} & \multicolumn{1}{l}{$p=-2$} & \multicolumn{1}{l}{$p=2$} & $p=2.5$ \\ \cline{3-10} 
   			&                                & \multicolumn{1}{l}{$\delta_{1RMLI}$}       & \multicolumn{1}{l}{$\delta_{1RMLI}$}     & \multicolumn{1}{l}{$\delta_{1RMLI}$}    & $\delta_{1RMLI}$      & \multicolumn{1}{l}{$\delta_{1RMLI}$}       & \multicolumn{1}{l}{$\delta_{1RMLI}$}     & \multicolumn{1}{l}{$\delta_{1RMLI}$}    & $\delta_{1RMLI}$      \\ \hline \hline
   			\multirow{7}{*}{(8,8)}       & 0.1                            & \multicolumn{1}{l}{50.92}    & \multicolumn{1}{l}{51.60}  & \multicolumn{1}{l}{57.81} & 58.70   & \multicolumn{1}{l}{55.59}    & \multicolumn{1}{l}{56.00}  & \multicolumn{1}{l}{59.65} & 56.00   \\ 
   			& 0.2                            & \multicolumn{1}{l}{45.81}    & \multicolumn{1}{l}{46.56}  & \multicolumn{1}{l}{53.69} & 54.75   & \multicolumn{1}{l}{52.98}    & \multicolumn{1}{l}{53.50}  & \multicolumn{1}{l}{58.03} & 53.50   \\ 
   			& 0.4                            & \multicolumn{1}{l}{40.62}    & \multicolumn{1}{l}{41.28}  & \multicolumn{1}{l}{48.08} & 49.17   & \multicolumn{1}{l}{44.99}    & \multicolumn{1}{l}{45.59}  & \multicolumn{1}{l}{50.97} & 45.59   \\ 
   			& 0.6                            & \multicolumn{1}{l}{39.14}    & \multicolumn{1}{l}{39.73}  & \multicolumn{1}{l}{45.82} & 46.82   & \multicolumn{1}{l}{39.38}    & \multicolumn{1}{l}{39.96}  & \multicolumn{1}{l}{45.39} & 39.96   \\ 
   			& 0.9                            & \multicolumn{1}{l}{38.74}    & \multicolumn{1}{l}{39.30}  & \multicolumn{1}{l}{44.86} & 45.76   & \multicolumn{1}{l}{34.71}    & \multicolumn{1}{l}{35.22}  & \multicolumn{1}{l}{40.25} & 35.22   \\ 
   			& 1.2                            & \multicolumn{1}{l}{38.71}    & \multicolumn{1}{l}{39.26}  & \multicolumn{1}{l}{44.69} & 45.54   & \multicolumn{1}{l}{32.36}    & \multicolumn{1}{l}{32.80}  & \multicolumn{1}{l}{37.30} & 32.80   \\ 
   			& 1.4                            & \multicolumn{1}{l}{38.70}    & \multicolumn{1}{l}{39.25}  & \multicolumn{1}{l}{44.65} & 45.49   & \multicolumn{1}{l}{31.48}    & \multicolumn{1}{l}{31.89}  & \multicolumn{1}{l}{36.07} & 31.89   \\ \hline \hline
   			\multirow{7}{*}{(9,10)}      & 0.1                            & \multicolumn{1}{l}{50.91}    & \multicolumn{1}{l}{51.54}  & \multicolumn{1}{l}{57.32} & 58.14   & \multicolumn{1}{l}{56.50}    & \multicolumn{1}{l}{56.88}  & \multicolumn{1}{l}{60.12} & 60.56   \\ 
   			& 0.2                            & \multicolumn{1}{l}{44.79}    & \multicolumn{1}{l}{45.45}  & \multicolumn{1}{l}{51.79} & 52.74   & \multicolumn{1}{l}{52.14}    & \multicolumn{1}{l}{52.61}  & \multicolumn{1}{l}{56.71} & 57.26   \\ 
   			& 0.4                            & \multicolumn{1}{l}{39.69}    & \multicolumn{1}{l}{40.26}  & \multicolumn{1}{l}{45.92} & 46.81   & \multicolumn{1}{l}{42.69}    & \multicolumn{1}{l}{43.20}  & \multicolumn{1}{l}{47.86} & 48.51   \\ 
   			& 0.6                            & \multicolumn{1}{l}{38.57}    & \multicolumn{1}{l}{39.07}  & \multicolumn{1}{l}{44.09} & 44.88   & \multicolumn{1}{l}{36.98}    & \multicolumn{1}{l}{37.46}  & \multicolumn{1}{l}{42.00} & 42.66   \\ 
   			& 0.9                            & \multicolumn{1}{l}{38.36}    & \multicolumn{1}{l}{38.84}  & \multicolumn{1}{l}{43.57} & 44.31   & \multicolumn{1}{l}{32.62}    & \multicolumn{1}{l}{33.03}  & \multicolumn{1}{l}{37.07} & 37.69   \\ 
   			& 1.2                            & \multicolumn{1}{l}{38.34}    & \multicolumn{1}{l}{38.83}  & \multicolumn{1}{l}{43.50} & 44.22   & \multicolumn{1}{l}{30.62}    & \multicolumn{1}{l}{30.97}  & \multicolumn{1}{l}{34.49} & 35.05   \\ 
   			& 1.4                            & \multicolumn{1}{l}{38.34}    & \multicolumn{1}{l}{38.83}  & \multicolumn{1}{l}{43.50} & 44.22   & \multicolumn{1}{l}{29.95}    & \multicolumn{1}{l}{30.27}  & \multicolumn{1}{l}{33.50} & 34.02   \\ \hline \hline
   			\multirow{7}{*}{(12,8)}      & 0.1                            & \multicolumn{1}{l}{50.21}    & \multicolumn{1}{l}{50.76}  & \multicolumn{1}{l}{55.72} & 56.42   & \multicolumn{1}{l}{55.26}    & \multicolumn{1}{l}{55.62}  & \multicolumn{1}{l}{58.76} & 59.19   \\ 
   			& 0.2                            & \multicolumn{1}{l}{45.87}    & \multicolumn{1}{l}{46.43}  & \multicolumn{1}{l}{51.70} & 52.47   & \multicolumn{1}{l}{52.08}    & \multicolumn{1}{l}{52.52}  & \multicolumn{1}{l}{56.39} & 56.91   \\ 
   			& 0.4                            & \multicolumn{1}{l}{42.99}    & \multicolumn{1}{l}{43.48}  & \multicolumn{1}{l}{48.22} & 48.94   & \multicolumn{1}{l}{44.70}    & \multicolumn{1}{l}{45.18}  & \multicolumn{1}{l}{49.49} & 50.09   \\ 
   			& 0.6                            & \multicolumn{1}{l}{42.54}    & \multicolumn{1}{l}{43.01}  & \multicolumn{1}{l}{47.45} & 48.12   & \multicolumn{1}{l}{40.35}    & \multicolumn{1}{l}{40.80}  & \multicolumn{1}{l}{44.95} & 45.55   \\ 
   			& 0.9                            & \multicolumn{1}{l}{42.51}    & \multicolumn{1}{l}{42.98}  & \multicolumn{1}{l}{47.38} & 48.04   & \multicolumn{1}{l}{37.38}    & \multicolumn{1}{l}{37.77}  & \multicolumn{1}{l}{41.46} & 42.01   \\ 
   			& 1.2                            & \multicolumn{1}{l}{42.51}    & \multicolumn{1}{l}{42.98}  & \multicolumn{1}{l}{47.38} & 48.04   & \multicolumn{1}{l}{36.24}    & \multicolumn{1}{l}{36.59}  & \multicolumn{1}{l}{39.88} & 40.37   \\ 
   			& 1.4                            & \multicolumn{1}{l}{42.51}    & \multicolumn{1}{l}{42.98}  & \multicolumn{1}{l}{47.38} & 48.04   & \multicolumn{1}{l}{35.91}    & \multicolumn{1}{l}{36.24}  & \multicolumn{1}{l}{39.34} & 39.80   \\ \hline \hline
   			\multirow{7}{*}{(14,15)}     & 0.1                            & \multicolumn{1}{l}{48.48}    & \multicolumn{1}{l}{48.93}  & \multicolumn{1}{l}{52.88} & 53.42   & \multicolumn{1}{l}{55.42}    & \multicolumn{1}{l}{55.71}  & \multicolumn{1}{l}{58.12} & 58.44   \\ 
   			& 0.2                            & \multicolumn{1}{l}{42.49}    & \multicolumn{1}{l}{42.91}  & \multicolumn{1}{l}{46.77} & 47.32   & \multicolumn{1}{l}{47.80}    & \multicolumn{1}{l}{48.15}  & \multicolumn{1}{l}{51.09} & 51.48   \\ 
   			& 0.4                            & \multicolumn{1}{l}{39.87}    & \multicolumn{1}{l}{40.22}  & \multicolumn{1}{l}{43.45} & 43.91   & \multicolumn{1}{l}{37.87}    & \multicolumn{1}{l}{38.20}  & \multicolumn{1}{l}{41.14} & 41.54   \\ 
   			& 0.6                            & \multicolumn{1}{l}{39.61}    & \multicolumn{1}{l}{39.94}  & \multicolumn{1}{l}{42.94} & 43.37   & \multicolumn{1}{l}{33.51}    & \multicolumn{1}{l}{33.79}  & \multicolumn{1}{l}{36.39} & 36.76   \\ 
   			& 0.9                            & \multicolumn{1}{l}{39.58}    & \multicolumn{1}{l}{39.91}  & \multicolumn{1}{l}{42.86} & 43.27   & \multicolumn{1}{l}{31.11}    & \multicolumn{1}{l}{31.34}  & \multicolumn{1}{l}{33.46} & 33.76   \\ 
   			& 1.2                            & \multicolumn{1}{l}{39.58}    & \multicolumn{1}{l}{39.91}  & \multicolumn{1}{l}{42.86} & 43.27   & \multicolumn{1}{l}{30.39}    & \multicolumn{1}{l}{30.59}  & \multicolumn{1}{l}{32.42} & 32.68   \\ 
   			& 1.4                            & \multicolumn{1}{l}{39.58}    & \multicolumn{1}{l}{39.91}  & \multicolumn{1}{l}{42.86} & 43.27   & \multicolumn{1}{l}{30.22}    & \multicolumn{1}{l}{30.41}  & \multicolumn{1}{l}{32.12} & 32.36   \\ \hline \hline
   			\multirow{7}{*}{(16,13)}     & 0.1                            & \multicolumn{1}{l}{48.70}    & \multicolumn{1}{l}{49.13}  & \multicolumn{1}{l}{52.84} & 53.35   & \multicolumn{1}{l}{55.23}    & \multicolumn{1}{l}{55.51}  & \multicolumn{1}{l}{57.91} & 58.23   \\ 
   			& 0.2                            & \multicolumn{1}{l}{43.94}    & \multicolumn{1}{l}{44.34}  & \multicolumn{1}{l}{47.93} & 48.44   & \multicolumn{1}{l}{48.70}    & \multicolumn{1}{l}{49.04}  & \multicolumn{1}{l}{51.95} & 52.34   \\ 
   			& 0.4                            & \multicolumn{1}{l}{42.23}    & \multicolumn{1}{l}{42.57}  & \multicolumn{1}{l}{45.70} & 46.14   & \multicolumn{1}{l}{39.97}    & \multicolumn{1}{l}{40.30}  & \multicolumn{1}{l}{43.19} & 43.59   \\ 
   			& 0.6                            & \multicolumn{1}{l}{42.08}    & \multicolumn{1}{l}{42.42}  & \multicolumn{1}{l}{45.42} & 45.84   & \multicolumn{1}{l}{36.28}    & \multicolumn{1}{l}{36.56}  & \multicolumn{1}{l}{39.12} & 39.48   \\ 
   			& 0.9                            & \multicolumn{1}{l}{42.08}    & \multicolumn{1}{l}{42.42}  & \multicolumn{1}{l}{45.40} & 45.82   & \multicolumn{1}{l}{34.46}    & \multicolumn{1}{l}{34.70}  & \multicolumn{1}{l}{36.85} & 37.15   \\ 
   			& 1.2                            & \multicolumn{1}{l}{42.08}    & \multicolumn{1}{l}{42.42}  & \multicolumn{1}{l}{45.40} & 45.82   & \multicolumn{1}{l}{34.01}    & \multicolumn{1}{l}{34.23}  & \multicolumn{1}{l}{36.16} & 36.43   \\ 
   			& 1.4                            & \multicolumn{1}{l}{42.08}    & \multicolumn{1}{l}{42.42}  & \multicolumn{1}{l}{45.40} & 45.82   & \multicolumn{1}{l}{33.91}    & \multicolumn{1}{l}{34.13}  & \multicolumn{1}{l}{35.99} & 36.25   \\ \hline \hline
   		\end{tabular}}
   	\end{table}
   	\begin{table}[H]
   		\caption{PRIs of the estimator $\delta_{2RMLI}$ with respect to RMLE when $\sigma_i$'s are known. }
   		\vspace{0.2cm}
   		\label{tab10.1}
   		\centering
   		\begin{tabular}{lllllllll}
   			\hline \hline
   			\multirow{3}{*}{$(n_1,n_2)$} & \multicolumn{4}{l}{$(\sigma_1,\sigma_2)=(1,1.5)$}                                                  & \multicolumn{4}{l}{$(\sigma_1,\sigma_2)=(3,2)$}                                                    \\ \cline{2-9} 
   			& \multicolumn{1}{l}{$p=-1$} & \multicolumn{1}{l}{$p=-0.5$} & \multicolumn{1}{l}{$p=0.5$} & $p=1$ & \multicolumn{1}{l}{$p=-1$} & \multicolumn{1}{l}{$p=-0.5$} & \multicolumn{1}{l}{$p=0.5$} & $p=1$ \\ \cline{2-9} 
   			& \multicolumn{1}{l}{$\delta_{2RMLI}$}     & \multicolumn{1}{l}{$\delta_{2RMLI}$}       & \multicolumn{1}{l}{$\delta_{2RMLI}$}      & $\delta_{2RMLI}$    & \multicolumn{1}{l}{$\delta_{2RMLI}$}     & \multicolumn{1}{l}{$\delta_{2RMLI}$}       & \multicolumn{1}{l}{$\delta_{2RMLI}$}      & $\delta_{2RMLI}$    \\ \hline \hline
   			(5,5)                        & \multicolumn{1}{l}{30.84}  & \multicolumn{1}{l}{31.37}    & \multicolumn{1}{l}{32.51}   & 33.10 & \multicolumn{1}{l}{39.93}  & \multicolumn{1}{l}{40.89}    & \multicolumn{1}{l}{43.00}   & 46.34 \\ 
   			(5,7)                        & \multicolumn{1}{l}{35.23}  & \multicolumn{1}{l}{35.76}    & \multicolumn{1}{l}{36.88}   & 37.47 & \multicolumn{1}{l}{42.74}  & \multicolumn{1}{l}{43.56}    & \multicolumn{1}{l}{45.35}   & 46.34 \\ 
   			(7,6)                        & \multicolumn{1}{l}{29.07}  & \multicolumn{1}{l}{29.47}    & \multicolumn{1}{l}{30.32}   & 30.78 & \multicolumn{1}{l}{39.02}  & \multicolumn{1}{l}{39.79}    & \multicolumn{1}{l}{41.45}   & 42.37 \\ 
   			(8,10)                       & \multicolumn{1}{l}{30.78}  & \multicolumn{1}{l}{34.35}    & \multicolumn{1}{l}{35.05}   & 35.42 & \multicolumn{1}{l}{42.34}  & \multicolumn{1}{l}{42.90}    & \multicolumn{1}{l}{44.08}   & 44.71 \\ 
   			(15,12)                      & \multicolumn{1}{l}{28.23}  & \multicolumn{1}{l}{28.41}    & \multicolumn{1}{l}{28.79}   & 28.98 & \multicolumn{1}{l}{38.73}  & \multicolumn{1}{l}{39.11}    & \multicolumn{1}{l}{39.88}   & 40.29 \\ \hline \hline
   		\end{tabular}
   	\end{table}
   	\begin{table}[H]
   		\caption{PRIs of the estimator $\delta_{2RMLI}$ with respect to RMLE when $\sigma_i$'s are known }
   		\label{tab10}
   		\centering
   		\begin{tabular}{lllllllll}
   			\hline \hline
   			\multirow{3}{*}{$(n_1,n_2)$} & \multicolumn{4}{l}{$(\sigma_1,\sigma_2)=(1,1.5)$}                                                  & \multicolumn{4}{l}{$(\sigma_1,\sigma_2)=(3,2)$}                                                    \\ \cline{2-9} 
   			& \multicolumn{1}{l}{$p=-2.5$} & \multicolumn{1}{l}{$p=-2$} & \multicolumn{1}{l}{$p=2$} & $p=2.5$ & \multicolumn{1}{l}{$p=-2.5$} & \multicolumn{1}{l}{$p=-2$} & \multicolumn{1}{l}{$p=2$} & $p=2.5$ \\ \cline{2-9} 
   			& \multicolumn{1}{l}{$\delta_{2RMLI}$}       & \multicolumn{1}{l}{$\delta_{2RMLI}$}     & \multicolumn{1}{l}{$\delta_{2RMLI}$}    & $\delta_{2RMLI}$      & \multicolumn{1}{l}{$\delta_{2RMLI}$}       & \multicolumn{1}{l}{$\delta_{2RMLI}$}     & \multicolumn{1}{l}{$\delta_{2RMLI}$}    & $\delta_{2RMLI}$      \\ \hline \hline
   			(8,8)                        & \multicolumn{1}{l}{30.17}    & \multicolumn{1}{l}{30.50}  & \multicolumn{1}{l}{33.48} & 33.91   & \multicolumn{1}{l}{38.86}    & \multicolumn{1}{l}{39.42}  & \multicolumn{1}{l}{44.92} & 45.77   \\ 
   			(9,10)                       & \multicolumn{1}{l}{33.91}    & \multicolumn{1}{l}{32.38}  & \multicolumn{1}{l}{35.13} & 35.54   & \multicolumn{1}{l}{40.52}    & \multicolumn{1}{l}{41.01}  & \multicolumn{1}{l}{45.77} & 46.50   \\ 
   			(12,8)                       & \multicolumn{1}{l}{24.96}    & \multicolumn{1}{l}{25.18}  & \multicolumn{1}{l}{27.07} & 27.32   & \multicolumn{1}{l}{35.19}    & \multicolumn{1}{l}{35.63}  & \multicolumn{1}{l}{39.83} & 40.44   \\ 
   			(14,15)                      & \multicolumn{1}{l}{32.05}    & \multicolumn{1}{l}{32.24}  & \multicolumn{1}{l}{33.97} & 34.21   & \multicolumn{1}{l}{40.95}    & \multicolumn{1}{l}{41.29}  & \multicolumn{1}{l}{44.34} & 44.77   \\ 
   			(16,13)                      & \multicolumn{1}{l}{28.22}    & \multicolumn{1}{l}{28.39}  & \multicolumn{1}{l}{29.86} & 30.06   & \multicolumn{1}{l}{38.26}    & \multicolumn{1}{l}{38.59}  & \multicolumn{1}{l}{41.56} & 41.56   \\ \hline \hline
   		\end{tabular}
   	\end{table}
   	
\subsection{ Improved estimation of $\mu_i$ when $\sigma_i$'s are equal but unknown} \label{sec3.2}
In this section, we consider the improved estimation of $\mu_i$ when $\sigma_1=\dots=\sigma_k=\sigma$ and $\sigma$ is unknown. The estimation problem will be studied with respect to the loss function
\begin{eqnarray}\label{lnloss2}
	L(\mu_{i},\sigma,\delta_{i})=\exp\left\{
	p\left(\frac{\delta_{i}-\mu_{i}}{\sigma}\right)\right\}
	-p\left(\frac{\delta_{i}-\mu_{i}}{\sigma}\right)-1,~p\neq0.
\end{eqnarray}
The joint complete and sufficient statistic for this estimation problem  is $(X_{i(1)},\dots,X_{k(1)}, T)$, where $T=\sum\limits_{i=1}^k T_i$. Also $X_{i(1)},\dots,X_{k(1)}, T$ are independently distributed. Further, $T\sim Gamma(n-k,\sigma)$.
Consider a group of transformations as $G_{a,b_i}=\{g_{a,b_i}(x)=ax_{ij}+b_i, j=1,\dots,n_i, i=1,\dots,k\}.$ After some simplification, we get the form of the $G_{a,b_i}$-equivariant estimators of $\mu_i$ as 
$$\delta_{\beta_i}(X_{i(1)},T)=X_{i(1)}+\beta_i T,$$
where $\beta_i$ is a real constant. The risk of the estimator $\delta_{\beta_i}(X_{i(1)},T)$ under the loss function (\ref{lnloss2})  is 
\begin{equation}
	R(\mu_i,\sigma,\delta_{\beta_i})=E\left[\exp\left\{p\left( \frac{X_{i(1)}+\beta_iT-\mu_i}{\sigma}\right) \right\}-p\left( {\frac{X_{i(1)}+\beta_iT-\mu_i}{\sigma}}\right) -1\right].
\end{equation}
To minimize the risk, differentiating $	R(\mu_i,\sigma,\delta_{\beta_i})$ with respect to $\beta_i,$ and then setting it to zero, we obtain the optimal choice of $\beta_i$ as
$$ \beta_{0i}=\frac{1}{p}\left( 1-\left(\frac{n_i}{n_i-p}\right) ^{\frac{1}{n-k+1}}\right),$$
provided $n_i>p.$ Hence, we have following theorem.
\begin{theorem}
	Under the Linex loss function in (\ref{lnloss2}), the best affine equivariant estimator of $\mu_{i}$ is $\delta_{\beta_{0i}}(X_{i(1)},T)$, where $ \beta_{0i}=\frac{1}{p}\left( 1-\left(\frac{n_i}{n_i-p}\right) ^{\frac{1}{n-k+1}}\right) .$
\end{theorem}

Now, our aim is to find an estimator that dominates $\delta_{\beta_{0i}}$. For this purpose, consider a class of estimators as 
 $$\delta_{\varphi_i}(\underline{X},T)=X_{i(1)}+\varphi_i(\underline{Y},T).$$
The risk of this estimator can be written as 
\begin{eqnarray*}
	R(\delta_{\varphi_{i}})&=&E\left[\exp\left\{p\left( \frac{X_{i(1)}+\varphi_i(\underline{Y},T)-\mu_i}{\sigma}\right)\right\} -p\left( \frac{X_{i(1)}+\varphi_i(\underline{Y},T)-\mu_i}{\sigma}\right)
	 -1\right]\\
	 &=&E^{\underline{Y},T}R^*(\delta_{\varphi_{i}}),
\end{eqnarray*}
where 
$$R^*(\delta_{\varphi_{i}})=E\left[\exp\left\{p\left( \frac{X_{i(1)}+\varphi_i(\underline{Y},T)-\mu_i}{\sigma}\right)\right\} -p\left( \frac{X_{i(1)}+\varphi_i(\underline{Y},T)-\mu_i}{\sigma}\right) -1 \big| T=t,\underline{Y}=\underline{y} \right].$$
After differentiating with respect to $\varphi_i(\underline{y},t)$ and equating it to zero we obtain 
\begin{equation}\label{eq3.14}
	\varphi_i(\underline{y},t)=\mu_{i}-\mu_0-\frac{\sigma}{p}\ln\left(\frac{n}{n-p}\right),
\end{equation}
provided $n>p$ and $\mu_0= \max\left\{ \mu_{i}, \max\limits_{1 \le j \le k, i\neq j}(\mu_{j}-y_{j})\right\}$. The supremum and infimum of  (\ref{eq3.14}) under the restriction $\mu_{1} \le\dots \le \mu_{k}$ are obtained as 
\begin{eqnarray}
\varphi_i^{*}(\underline{y})=\sup\limits_{\mu_\le\dots\le\mu_k}\ {\varphi_i(\underline{y},t)}&=&
	\left\{\begin{array}{ll}
		\displaystyle \min(0,y_{i+1},\ldots,y_{k}),
		& \textrm{if $i=1,\ldots,k-1$}\\
		\displaystyle 0  &
		\textrm{if $i=k$}
	\end{array} \right. \nonumber 
\end{eqnarray}
and $$\inf\limits_{\mu_\le\dots\le\mu_k} \ {\varphi_i(\underline{y},t)} = -\infty.$$
 Now, define 
\begin{eqnarray}
 {\varphi_{0i}(\underline{y},t)}&=&
	\left\{\begin{array}{ll}
		\displaystyle  \varphi ^*(\underline{y},t),
		& \textrm{if $\varphi_i(\underline{y},t)\ge \varphi ^*(\underline{y},t)$ }\\
		\displaystyle \varphi_i(\underline{y},t)  &
		\textrm{otherwise.}
	\end{array} \right.
\end{eqnarray}

\begin{theorem}
	Suppose $\sigma_{1}=\dots=\sigma_k = \sigma$ and $\mu_1,\dots,\mu_k$ satisfy the order restriction $\mu_{1} \le \dots \le \mu_{k} $. Then, for estimating $\mu_i$, the risk of the estimator $X_{i(1)} + \min\left\lbrace \varphi_i^{*}(\underline{Y},T),\varphi_i(\underline{Y},T)\right\rbrace $ is nowhere larger than $X_{i(1)} + \varphi_i(\underline{Y},T)$ under the Linex loss function provided $P\left(\varphi_i^{*}(\underline{Y},T)<\varphi_i(\underline{Y},T)\right)\ne0$.
\end{theorem}

As an application to the above theorem, we derive the  following corollary that provides an estimator which dominates the BAEE. 
\begin{corollary}
Let $\mu_{1} \le \dots \le \mu_{k}$ and $\sigma_i =\sigma$ for $i =1,\dots, k$. For  estimating $\mu_{i}$, $i =1,\dots,k$ the estimator 
\begin{eqnarray}
	{\delta_{\varphi_{0i}}(\underline{Y},T)}&=&
	\left\{\begin{array}{ll}
		\displaystyle X_{i(1)} + \min\left\lbrace \beta_{0i}T,\min(0,y_{i+1},\ldots,y_{k} \right\rbrace  ,
		& \textrm{for $i =1 \dots k-1$ }\\
		\displaystyle X_{k(1)}+\min\{\beta_{0k}T, 0\},  &
		\textrm{for $i =k$}
	\end{array} \right.
\end{eqnarray}
dominates the best location equivariant estimator $\delta_{\beta_{0i}}$ under the Linex loss function.

\end{corollary}
\begin{remark}
	For estimating $\mu_2$, it is observed that $\delta_{\varphi_{0i}}$ improves upon the estimator $\delta_{2RML}$.
\end{remark}
\subsection{A simulation study}
In Table \ref{tab8}, we have tabulated the PRIs  of  $(\delta_{\beta_{01}},\delta_{\beta_{02}})$ with respect to MLE. In Table \ref{tab9}, we have provided the PRI  of $\delta_{\varphi_{01}}$ with respect to the BAEE. In this case, $\delta_{\varphi_{02}}$ coincides with the best affine equivariant estimator. The risk functions of  $\delta_{\beta_{01}}$ and $\delta_{\beta_{02}}$ do not depend on the unknown parameter.  We have the following observations from Tables \ref{tab8} and \ref{tab9}.
\begin{itemize}
	\item [(i)]  PRIs  of  $\delta_{\beta_{01}}$ and $\delta_{\beta_{02}}$  increase as $p$ increases.  As the value of $p$ becomes large the estimators  $\delta_{\beta_{01}}$ and $\delta_{\beta_{02}}$ perform significantly better than MLE. 
	\item[(ii)]	The PRI value of $\delta_{\varphi_{01}}$ increases as $p$ increases. For large values of $p$ with small values  $\mu_{1}-\mu_{2}$, the estimator  $\delta_{\varphi_{01}}$ provides significantly better improvement over $\delta_{\beta_{01}}$. Further, as $\mu_2-\mu_1$ and $\sigma$ increase, PRI of $\delta_{\varphi_{01}}$ decreases. 

\end{itemize}\begin{table}[H]
\caption{PRIs of the BAEEs $\delta_{\beta_{01}}$ and $\delta_{\beta_{02}}$ with respect to the MLE, when $\sigma_{i}$'s  are unknown but equal.}
\label{tab8}
\centering
\resizebox{18cm}{1.5cm}{
	
	\begin{tabular}{lllllllll}
		\hline \hline
		\multirow{2}{*}{$(n_1,n_2)$} & $p=-4$                      & $p=-2$                      & $p=-1$                      & $p=-0.5$                    & $p=0.5$                     & $p=1$                       & $p=2$                       & $p=4$                       \\ \cline{2-9} 
		& $(\delta_{\beta_{01}},\delta_{\beta_{02}})$ & $(\delta_{\beta_{01}},\delta_{\beta_{02}})$ & $(\delta_{\beta_{01}},\delta_{\beta_{02}})$ & $(\delta_{\beta_{01}},\delta_{\beta_{02}})$ & $(\delta_{\beta_{01}},\delta_{\beta_{02}})$ & $(\delta_{\beta_{01}},\delta_{\beta_{02}})$ & $(\delta_{\beta_{01}},\delta_{\beta_{02}})$ & $(\delta_{\beta_{01}},\delta_{\beta_{02}})$ \\ \hline\hline
		$(5,5)$                      & (35.15, 34.95)              & (39.11, 38.85)              & (41.61, 41.31)              & (43.04, 42.71)              & (46.36, 45.94)              & (48.30, 47.81)              & (52.93, 52.10)              & (66.12, 65.38)              \\ 
		$(5,7)$                      & (35.78, 37.81)              & (39.74, 41.04)              & (42.23, 42.97)              & (43.66, 44.04)              & (47.01, 46.42)              & (49.00, 47.76)              & (53.97, 50.80)              & (70.34, 59.02)              \\ 
		$(7,6)$                      & (38.49, 37.70)              & (41.73, 41.28)              & (43.66, 43.48)              & (44.74, 44.71)              & (47.13, 47.50)              & (48.49, 49.12)              & (51.59, 52.92)              & (60.17, 64.25)              \\ 
		$(8,10)$                     & (40.10, 41.17)              & (43.08, 43.70)              & (44.82, 45.14)              & (45.78, 45.92)              & (47.87, 47.60)              & (49.03, 48.51)              & (51.64, 50.50)              & (58.57, 55.33)              \\ 
		$(15,12)$                    & (44.01, 43.16)              & (45.86, 45.35)              & (46.87, 46.58)              & (47.40, 47.23)              & (48.53, 48.61)              & (49.13, 49.34)              & (50.39, 50.92)              & (53.26, 54.59)              \\ \hline \hline
\end{tabular}}
\end{table}

\newpage
\begin{table}[H]
			\caption{PRIs of the estimator $\delta_{\varphi_{01}}$ with respect to the BAEE when $\sigma_{i}$'s  are unknown but equal. }
			\label{tab9}
	\centering
	\resizebox{16cm}{10.5cm}{

		\begin{tabular}{lllllllllll}
			\hline \hline
			\multirow{2}{*}{$(n_1,n_2)$} & \multirow{2}{*}{$\mu_2-\mu_1$} & \multirow{2}{*}{$\sigma$} & $p=-4$                  & $p=-2$                  & $p=-1$                  & $p=-0.5$                & $p=0.5$                 & $p=1$                   & $p=2$                   & $p=4$                   \\ \cline{4-11} 
			&                                &                           & $\delta_{\varphi_{01}}$ & $\delta_{\varphi_{01}}$ & $\delta_{\varphi_{01}}$ & $\delta_{\varphi_{01}}$ & $\delta_{\varphi_{01}}$ & $\delta_{\varphi_{01}}$ & $\delta_{\varphi_{01}}$ & $\delta_{\varphi_{01}}$ \\ \hline
			\multirow{7}{*}{$(5,5)$}     & 0.1                            & 0.7                       & 35.68                   & 38.38                   & 40.07                   & 41.01                   & 43.12                   & 44.29                   & 44.29                   & 51.37                   \\ 
			& 0.2                            & 1                         & 28.46                   & 31.42                   & 33.36                   & 34.47                   & 37.03                   & 38.49                   & 38.49                   & 49.48                   \\ 
			& 0.4                            & 1.5                       & 21.45                   & 24.39                   & 26.42                   & 27.62                   & 30.51                   & 32.25                   & 32.25                   & 47.22                   \\ 
			& 0.6                            & 2                         & 18.53                   & 21.37                   & 23.39                   & 24.61                   & 27.57                   & 29.39                   & 29.39                   & 46.06                   \\ 
			& 0.9                            & 2.3                       & 12.23                   & 14.63                   & 16.49                   & 17.66                   & 20.65                   & 22.56                   & 22.56                   & 42.94                   \\ 
			& 1.2                            & 2.8                       & 10.28                   & 12.50                   & 14.25                   & 15.36                   & 18.25                   & 20.14                   & 20.14                   & 41.74                   \\ 
			& 1.4                            & 3                         & 8.61                    & 10.60                   & 12.23                   & 13.28                   & 16.06                   & 17.91                   & 17.91                   & 40.54                   \\ \hline \hline
			\multirow{7}{*}{$(5,7)$}     & 0.1                            & 0.7                       & 41.10                   & 43.81                   & 45.43                   & 46.31                   & 48.18                   & 48.18                   & 50.71                   & 42.42                   \\ 
			& 0.2                            & 1                         & 32.98                   & 36.13                   & 38.13                   & 39.26                   & 41.74                   & 41.74                   & 45.60                   & 40.14                   \\ 
			& 0.4                            & 1.5                       & 25.07                   & 28.31                   & 30.51                   & 31.78                   & 34.72                   & 34.72                   & 39.78                   & 37.31                   \\ 
			& 0.6                            & 2                         & 21.73                   & 24.92                   & 27.14                   & 28.44                   & 31.50                   & 31.50                   & 36.99                   & 35.86                   \\ 
			& 0.9                            & 2.3                       & 14.52                   & 17.28                   & 19.33                   & 20.60                   & 23.72                   & 23.72                   & 29.92                   & 31.82                   \\ 
			& 1.2                            & 2.8                       & 12.24                   & 14.77                   & 16.72                   & 17.95                   & 21.01                   & 21.01                   & 27.32                   & 30.18                   \\ 
			& 1.4                            & 3                         & 10.25                   & 12.56                   & 14.38                   & 15.55                   & 18.51                   & 18.51                   & 24.81                   & 28.51                   \\ \hline \hline
			\multirow{7}{*}{$(7,6)$}     & 0.1                            & 0.7                       & 27.58                   & 29.86                   & 31.25                   & 32.02                   & 33.70                   & 34.62                   & 36.54                   & 39.82                   \\ 
			& 0.2                            & 1                         & 19.90                   & 22.19                   & 23.65                   & 24.47                   & 26.33                   & 27.38                   & 29.70                   & 34.40                   \\ 
			& 0.4                            & 1.5                       & 13.19                   & 15.21                   & 16.58                   & 17.37                   & 19.24                   & 20.32                   & 22.84                   & 28.69                   \\ 
			& 0.6                            & 2                         & 10.66                   & 12.50                   & 13.78                   & 14.54                   & 16.34                   & 17.41                   & 19.92                   & 26.08                   \\ 
			& 0.9                            & 2.3                       & 5.88                    & 7.18                    & 8.15                    & 8.75                    & 10.24                   & 11.16                   & 13.41                   & 19.65                   \\ 
			& 1.2                            & 2.8                       & 4.56                    & 5.66                    & 6.49                    & 7.02                    & 8.35                    & 9.18                    & 11.30                   & 17.41                   \\ 
			& 1.4                            & 3                         & 3.50                    & 4.41                    & 5.13                    & 5.59                    & 6.77                    & 7.53                    & 9.49                    & 15.34                   \\ \hline \hline
			\multirow{7}{*}{$(8,10)$}    & 0.1                            & 0.7                       & 30.20                   & 32.58                   & 32.58                   & 34.76                   & 34.76                   & 37.34                   & 39.29                   & 43.47                   \\ 
			& 0.2                            & 1                         & 21.02                   & 23.40                   & 23.40                   & 25.73                   & 25.73                   & 28.68                   & 31.05                   & 36.69                   \\ 
			& 0.4                            & 1.5                       & 13.31                   & 15.38                   & 15.38                   & 17.55                   & 17.55                   & 20.51                   & 23.05                   & 29.71                   \\ 
			& 0.6                            & 2                         & 10.53                   & 12.39                   & 12.39                   & 14.42                   & 14.42                   & 17.27                   & 19.77                   & 26.61                   \\ 
			& 0.9                            & 2.3                       & 5.46                    & 6.72                    & 6.72                    & 8.21                    & 8.21                    & 10.51                   & 12.71                   & 19.45                   \\ 
			& 1.2                            & 2.8                       & 4.14                    & 5.18                    & 5.18                    & 6.47                    & 6.47                    & 8.53                    & 10.56                   & 17.01                   \\ 
			& 1.4                            & 3                         & 3.12                    & 3.96                    & 3.96                    & 5.05                    & 5.05                    & 6.83                    & 8.66                    & 14.78                   \\ \hline \hline
			\multirow{7}{*}{$(15,12)$}   & 0.1                            & 0.7                       & 11.94                   & 13.07                   & 13.72                   & 14.07                   & 14.82                   & 15.23                   & 16.08                   & 17.97                   \\ 
			& 0.2                            & 1                         & 5.55                    & 6.28                    & 6.72                    & 6.96                    & 7.49                    & 7.78                    & 8.43                    & 9.96                    \\ 
			& 0.4                            & 1.5                       & 2.09                    & 2.45                    & 2.68                    & 2.81                    & 3.10                    & 3.27                    & 3.65                    & 4.60                    \\ 
			& 0.6                            & 2                         & 1.27                    & 1.52                    & 1.68                    & 1.78                    & 2.00                    & 2.12                    & 2.41                    & 3.16                    \\ 
			& 0.9                            & 2.3                       & 0.32                    & 0.39                    & 0.45                    & 0.48                    & 0.56                    & 0.60                    & 0.71                    & 1.01                    \\ 
			& 1.2                            & 2.8                       & 0.17                    & 0.22                    & 0.25                    & 0.27                    & 0.32                    & 0.34                    & 0.41                    & 0.61                    \\ 
			& 1.4                            & 3                         & 0.08                    & 0.10                    & 0.12                    & 0.13                    & 0.15                    & 0.16                    & 0.19                    & 0.27                    \\ \hline \hline
		\end{tabular}    }
\end{table}
\newpage 

\newpage
\subsection{Improved estimation $\mu_i$ when all $\sigma_{i}$'s are unknown and unequal} \label{sec3.3}
In this section, we study the estimation of $\mu_i$ under the order restriction $\mu_1\le \dots\le \mu_k$ when $\sigma_i$'s are unknown and unequal with respect to the loss function (\ref{linloss1}).  The joint complete and sufficient statistic for this estimation problem is $(X_{1(1)},\dots,X_{k(1)}, T_1,T_2\dots,T_k)$.
Consider an affine group of transformations as $$G_{a_i,b_i}=\{g_{a_i,b_i}(x)=ax_{ij}+b_i, j=1,\dots,n_i\},$$ where
$i=1,\dots,k$. The form of  affine equivariant estimators is obtained as
$$\delta_{\kappa_i}(X_{i(1)},T_i)=X_{i(1)}+\kappa_i T_i.$$
The risk of estimator under Linex loss function given in (\ref{linloss1}) is 
\begin{equation}\label{ss3e1}
	R(\mu_i,\sigma_i, \delta_{\kappa_i})=E\left[ \exp \left\{p\left(\frac{X_{i(1)}+\kappa_i T_i-\mu_{i}}{\sigma_{i}} \right)\right\} - p\left(\frac{X_{i(1)}+\kappa_i T_i-\mu_{i}}{\sigma_{i}} \right) -1 \right].
\end{equation}
After differentiating (\ref{ss3e1}) with respect to $\kappa_i$, and then equating  to zero, the value of $\kappa_i$ is obtained as
$\kappa_i=\frac{1}{p}\left(1-\left( {\frac{n_i}{n_i-p}}\right)^{\frac{1}{n_i}} \right),$ provided $n_i>p$.
Thus, we have  the following theorem.
\begin{theorem}
		Under the Linex loss function (\ref{linloss1}) the best affine equivariant estimator of $\mu_{i}$ is $\delta_{\kappa_{0i}}(X_{i(1)},T_i)$, where $\kappa_{0i}$ is 
	$\kappa_{0i}=\frac{1}{p}\left(1-\left( {\frac{n_i}{n_i-p}}\right)^{\frac{1}{n_i}} \right)$,  $n_i>p.$
\end{theorem}

Now our goal is to derive an estimator which  improves upon  $\delta_{\kappa_{0i}}$ with respect to the Linex loss function. For this purpose, we consider a class of estimators as 
$$\delta_{\zeta_i}= X_{i(1)}+\zeta_i(\underline{Y},W_i ).$$
The risk  function of $\delta_{\zeta_i}$ can be expressed as 
\begin{eqnarray*}
		R(\underline{\mu}, \sigma_i,\delta_{\zeta_i})&=&E\left[ \exp\left\lbrace p\left( \frac{X_{i(1)}+\zeta_i(\underline{Y},W_i)-\mu_{i}}{\sigma_{i}}\right)\right\rbrace  -p\left( \frac{X_{i(1)}+\zeta_i(\underline{Y},W_i)-\mu_{i}}{\sigma_{i}}\right) -1 \right] \\
&=&ER^*(\underline{\mu}, \delta_{\zeta_i}),
\end{eqnarray*}
where $\underline{\mu}=(\mu_1,\dots,\mu_k)$ and 
\begin{equation*}
R^*(\underline{\mu}, \sigma_i, \delta_{\zeta_i})=E\left[ \exp\left\lbrace p\left( \frac{X_{i(1)}+\zeta_i(\underline{Y},W_i)-\mu_{i}}{\sigma_{i}}\right)\right\rbrace  -p\left( \frac{X_{i(1)}+\zeta_i(\underline{Y},W_i)-\mu_{i}}{\sigma_{i}}\right) -1\big|\underline{Y},W_i=w_i  \right].
\end{equation*}
It can be easily noticed that the conditional risk $R^*(\underline{\mu}, \sigma_i, \delta_{\zeta_i})$ is minimized at 
\begin{equation}\label{eq3.19}
	\zeta_i(\underline{y},w_i) =\mu_i -\mu_0-\frac{\sigma_{i}}{p}\ln \left(\frac{ n}{ n-p}\right),
\end{equation}
where $\mu_0= \max\left\lbrace \mu_{i}, {\max\limits_{1 \le j \le k, i\neq j}(\mu_{j}-y_{j})}\right\rbrace $. Under the order restriction $\mu_{1} \le \dots \le \mu_{k}$, the supremum and infimum for $\zeta_i(\underline{y},w_i)$ can be obtained as 

\begin{eqnarray}
	{\zeta_i}^{*}(\underline{y})=\sup	 \ {\zeta_i(\underline{y},w_i)}&=&
	\left\{\begin{array}{ll}
		\displaystyle \min(0,y_{i+1},\ldots,y_{k}),
		& \textrm{if $i=1,\ldots,k-1$}\\
		\displaystyle 0 &
		\textrm{if $i=k$}
	\end{array} \right. \nonumber 
\end{eqnarray}
and
$$\inf \ {\zeta_i(\underline{Y},W_i)} = -\infty.$$
Define 
\begin{eqnarray}
	 {\zeta_i}^{*}(\underline{y}),=\ {\zeta_{0i}(\underline{y},w_i)}&=&
	\left\{\begin{array}{ll}
		\displaystyle \zeta_{i}^*(\underline{y},w_i),
		& \textrm{if $\zeta_{i}(\underline{y},w_i)\ge\zeta_{i}^*(\underline{y},w_i)$}\\
		\displaystyle \zeta_{i}(\underline{y},w_i)  &
		\textrm{otherwise.}
	\end{array} \right. \nonumber
\end{eqnarray}
\begin{theorem} 
	Suppose $\sigma_{i}, i,=1,\dots,k $ are unknown and $\mu_1, \dots, \mu_{k}$ are associated to he order restriction $\mu_{1} \le \dots < \mu_{k}.$ Then the risk of the estimator $X_{i(1)}+ \min\left\lbrace {{\zeta}^{*}(\underline{Y},W_i),{\zeta}(\underline{Y},W_i)}\right\rbrace $ is nowhere larger than $X_{i(1)} + {{\zeta}(\underline{Y},W_i)}$ under the Linex loss (\ref{linloss1}),  provided 
	$P_{\underline{\mu}}(\zeta_{i0}(\underline{y},w_i)\neq\zeta_{i}(\underline{Y},W_i))>0$.
\end{theorem}
An application of the preceding theorem leads to the following corollary. 
\begin{corollary}
	Consider the order restriction $\mu_{1} \le \dots \le \mu_{k}$ and all $\sigma_i $ for $i =1, \dots, k$ are unknown. To  estimate $\mu_{i}$, $i =1, \dots, k$, the estimator 
	\begin{eqnarray}
		{\delta_{\zeta_{0i}}(\underline{Y},W_i)}&=&
		\left\{\begin{array}{ll}
			\displaystyle X_{i(1)}+\min\{\kappa_{0i}T_i, \min(0,Y_{i+1},\dots,Y_{k})\}  ,
			& \textrm{for $i =1,\dots, k-1$ }\\
			\displaystyle X_{k(1)}+\min\{\kappa_{0k}T_{k},0\}\   &
			\textrm{for $i=k$}
		\end{array} \right.
	\end{eqnarray}
	dominates the estimator $\delta_{\kappa_{0i}}$ under the Linex loss function (\ref{linloss1}).
\end{corollary}
\subsection{A simulation study}
In Table \ref{tab15} we have tabulated the PRI of $\delta_{\zeta_{01}}$ with respect to $\delta_{\kappa_{01}}$ by taking various samples sizes and $p$ values. Following observations have been drawn from Table \ref{tab15}.
\begin{itemize}
	\item [(i)] PRI of $\delta_{\zeta_{01}}$ increases as $p$ increases and decreases as $\sigma_{1}$ and $\sigma_{2}$ increase.
 \item[(ii)]PRI decreases as sample size and $\mu_{2}-\mu_{1}$ increase, moreover for closed values of $\mu_{1}$ and $\mu_{2}$, the estimator  $\delta_{\zeta_{01}}$ improves significantly better. 
\end{itemize}
\begin{table}[H]
	\caption{PRIs of $\delta_{\zeta_{01}}$ with respect to the BAEE when $\sigma_{i}$'s are unknown and unequal.}
	\vspace{0.2cm}
	\label{tab15}
	\centering
	\resizebox{16cm}{10.5cm}{
		\begin{tabular}{llllllllllll}
			\hline \hline
			\multirow{2}{*}{$(n_1,n_2)$} & \multirow{2}{*}{$\mu_2-\mu_1$} & \multirow{2}{*}{$\sigma_1$} & \multirow{2}{*}{$\sigma_2$} & $p=-4$                 & $p=-2$                 & $p=-1$                 & $p=-0.5$               & $p=0.5$                & $p=1$                  & $p=2$                  & $p=4$                  \\ \cline{5-12} 
			&                                &                             &                             & $\delta_{\zeta_{01}}$ & $\delta_{\zeta_{01}}$ & $\delta_{\zeta_{01}}$ & $\delta_{\zeta_{01}}$ & $\delta_{\zeta_{01}}$ & $\delta_{\zeta_{01}}$ & $\delta_{\zeta_{01}}$ & $\delta_{\zeta_{01}}$ \\ \hline \hline
											\multirow{7}{*}{$(5,5)$}     & 0.1                            & 0.7                         & 0.5                         & 42.24                  & 45.18                  & 46.97                  & 47.96                  & 50.11                  & 51.28                  & 53.83                  & 56.88                  \\ 
											& 0.2                            & 1                           & 0.9                         & 30.75                  & 34.01                  & 36.16                  & 37.39                  & 40.26                  & 41.92                  & 45.80                  & 53.58                  \\ 
											& 0.4                            & 1.5                         & 1.6                         & 21.45                  & 24.56                  & 26.74                  & 28.05                  & 31.21                  & 33.14                  & 37.93                  & 50.15                  \\ 
											& 0.6                            & 2                           & 1.9                         & 19.62                  & 22.73                  & 24.95                  & 26.29                  & 29.58                  & 31.60                  & 36.62                  & 49.66                  \\ 
											& 0.9                            & 2.3                         & 2.5                         & 12.16                  & 14.65                  & 16.58                  & 17.80                  & 20.96                  & 23.02                  & 28.54                  & 45.72                  \\ 
											& 1.2                            & 2.8                         & 3                           & 10.30                  & 12.58                  & 14.39                  & 15.56                  & 18.65                  & 20.70                  & 26.31                  & 44.56                  \\ 
											& 1.4                            & 3                           & 3.5                         & 8.22                   & 10.22                  & 11.86                  & 12.95                  & 15.88                  & 17.86                  & 23.44                  & 42.78                  \\ \hline \hline
											\multirow{7}{*}{$(5,7)$}     &  0.1                            & 0.7                         & 0.5                         & 47.51                  & 50.40                  & 52.07                  & 52.96                  & 54.80                  & 55.72                  & 57.25                  & 49.97                  \\ 
											& 0.2                            & 1                           & 0.9                         & 35.44                  & 38.87                  & 41.04                  & 42.25                  & 44.97                  & 46.44                  & 49.41                  & 46.39                  \\ 
											& 0.4                            & 1.5                         & 1.6                         & 25.27                  & 28.69                  & 31.02                  & 32.39                  & 35.59                  & 37.43                  & 41.44                  & 42.34                  \\ 
											& 0.6                            & 2                           & 1.9                         & 22.91                  & 26.32                  & 28.70                  & 30.10                  & 33.43                  & 35.37                  & 39.67                  & 41.52                  \\ 
											& 0.9                            & 2.3                         & 2.5                         & 14.46                  & 17.30                  & 19.44                  & 20.77                  & 24.09                  & 26.14                  & 30.98                  & 36.53                  \\ 
											& 1.2                            & 2.8                         & 3                           & 12.26                  & 14.87                  & 16.91                  & 18.18                  & 21.43                  & 23.45                  & 28.38                  & 34.94                  \\ 
											& 1.4                            & 3                           & 3.5                         & 9.89                   & 12.17                  & 14.01                  & 15.18                  & 18.24                  & 20.21                  & 25.16                  & 32.70                  \\ \hline \hline
											\multirow{7}{*}{$(7,6)$}     & 0.1                            & 0.7                         & 0.5                         & 32.85                  & 35.37                  & 36.88                  & 37.69                  & 39.43                  & 40.35                  & 42.25                  & 45.16                  \\ 
											& 0.2                            & 1                           & 0.9                         & 21.43                  & 23.90                  & 25.46                  & 26.34                  & 28.32                  & 29.42                  & 31.82                  & 36.72                  \\ 
											& 0.4                            & 1.5                         & 1.6                         & 13.03                  & 15.06                  & 16.44                  & 17.24                  & 19.14                  & 20.24                  & 22.81                  & 28.93                  \\ 
											& 0.6                            & 2                           & 1.9                         & 11.22                  & 13.16                  & 14.50                  & 15.30                  & 17.18                  & 18.28                  & 20.89                  & 27.31                  \\ 
											& 0.9                            & 2.3                         & 2.5                         & 5.71                   & 6.99                   & 7.94                   & 8.52                   & 9.98                   & 10.88                  & 13.13                  & 19.56                  \\ 
											& 1.2                            & 2.8                         & 3                           & 4.46                   & 5.54                   & 6.36                   & 6.88                   & 8.19                   & 9.02                   & 11.14                  & 17.37                  \\ 
											& 1.4                            & 3                           & 3.5                         & 3.26                   & 4.11                   & 4.79                   & 5.22                   & 6.32                   & 7.04                   & 8.88                   & 14.59                  \\ \hline \hline
											\multirow{7}{*}{$(8,10)$}    & 0.1                            & 0.7                         & 0.5                         & 34.69                  & 37.25                  & 38.74                  & 39.55                  & 41.28                  & 42.21                  & 44.18                  & 48.29                  \\ 
											& 0.2                            & 1                           & 0.9                         & 22.41                  & 24.97                  & 26.58                  & 27.49                  & 29.53                  & 30.67                  & 33.21                  & 39.29                  \\ 
											& 0.4                            & 1.5                         & 1.6                         & 13.31                  & 15.45                  & 16.90                  & 17.75                  & 19.72                  & 20.87                  & 23.57                  & 30.73                  \\ 
											& 0.6                            & 2                           & 1.9                         & 11.08                  & 13.08                  & 14.45                  & 15.26                  & 17.18                  & 18.31                  & 21.00                  & 28.35                  \\ 
											& 0.9                            & 2.3                         & 2.5                         & 5.48                   & 6.76                   & 7.72                   & 8.31                   & 9.77                   & 10.69                  & 12.99                  & 20.07                  \\ 
											& 1.2                            & 2.8                         & 3                           & 4.18                   & 5.25                   & 6.07                   & 6.58                   & 7.88                   & 8.70                   & 10.82                  & 17.64                  \\ 
											& 1.4                            & 3                           & 3.5                         & 3.02                   & 3.87                   & 4.54                   & 4.97                   & 6.07                   & 6.79                   & 8.68                   & 15.03                  \\ \hline \hline
											\multirow{7}{*}{$(15,12)$}   & 0.1                            & 0.7                         & 0.5                         & 14.26                  & 15.58                  & 16.34                  & 16.74                  & 17.61                  & 18.07                  & 19.04                  & 21.20                  \\ 
											& 0.2                            & 1                           & 0.9                         & 6.01                   & 6.80                   & 7.28                   & 7.55                   & 8.13                   & 8.45                   & 9.15                   & 10.80                  \\ 
											& 0.4                            & 1.5                         & 1.6                         & 2.07                   & 2.44                   & 2.67                   & 2.80                   & 3.10                   & 3.28                   & 3.66                   & 4.64                   \\ 
											& 0.6                            & 2                           & 1.9                         & 1.33                   & 1.59                   & 1.76                   & 1.86                   & 2.09                   & 2.22                   & 2.52                   & 3.31                   \\ 
											& 0.9                            & 2.3                         & 2.5                         & 0.31                   & 0.39                   & 0.44                   & 0.48                   & 0.55                   & 0.60                   & 0.70                   & 1.01                   \\ 
											& 1.2                            & 2.8                         & 3                           & 0.17                   & 0.21                   & 0.25                   & 0.27                   & 0.31                   & 0.34                   & 0.41                   & 0.60                   \\ 
											& 1.4                            & 3                           & 3.5                         & 0.08                   & 0.10                   & 0.11                   & 0.12                   & 0.14                   & 0.16                   & 0.19                   & 0.28                   \\ \hline \hline
										\end{tabular}
									}
								\end{table}

\section{Application to life testing sampling schemes}\label{sec4}
In the previous sections, we have derived the results of $\mu_{i}, i=1,\dots k$ by considering various scenarios under i.i.d. samples. In this section, we address the same estimation problem for three important life-testing schemes, such as $(i)$ Type-II censoring  $(ii)$ progressive Type-II censoring, and (iii) record values. 
\subsection{Type-II censoring}
In the survival analysis and reliability studies, it is not always possible to observe all test units until failure. Various situations may occur,  such as units being withdrawn, lost, or deliberately removed from the study before the experiment concludes.  Such incomplete observations lead to censored data. One of the most common censoring schemes is Type-II censoring. It occurs when the experimenter decides to terminate the experiment after a specific number of failures. For more details, one can refer to \cite{balakrishnan2000progressive}. Suppose a  random sample of size $n_i$ is drawn from $Exp(\mu_i,\sigma_{i})$ distribution and observations are in order $X_{i(1)} \le \dots \le X_{i(n_i)}$, for $i=1,\dots,k$. Under the Type-II censoring we consider only first $m_i$ ordered observations such as $X_{i(1)},\dots,  X_{i(m_i)}$, where $m_i<n_i$. From \cite{epstein1956simple} for $i=1,\dots,k$, we have  $X_{i(1)} \sim E( \mu_{i},\frac{\sigma_{i}}{n_i})$ and $T_i \sim \Gamma (m_i-1,\sigma_{i})$ with $T_i = \sum_{j=1}^{n_i} (X_{i(j)} - X_{i(1)}) + (n_i - m_i)(X_{i(n_i)} - X_{i(1)}) $. 
\begin{itemize}
	\item [(i)] For estimating $\mu_i$ with $\sigma_1\le \dots\le \sigma_k$ we have  $(X_{1(1)}, \dots, X_{k(1)}, T_1, \dots  T_k) $ is a joint complete sufficient statistic and it can be easily seen that this setup is similar to (\ref{m1}). As a consequence, improved estimators for $\mu_{i}$ can be derived similar to Section \ref{sec2}.

\item [(ii)] Now for estimating $\mu_i$ with the ordered restriction $\mu_1\le \dots \le \mu_k$ we have the following cases.  
$(a)$ When $\sigma_{i}$ are known $(X_{1(1)}, \dots, X_{k(1)})$ is complete and sufficient statistics, and this setup also aligns to Section \ref{sec3.1}. Therefore, improved estimators can be derived similar to Section \ref{sec3.1}.
$(b)$  When $\sigma_{i}$'s are equal but unknown we have joint sufficient statistic as $(X_{i(1)}, \dots X_{k(1)}, T)$ with $T=\sum_{i=1}^kT_i$.  This case  is similar to Section to \ref{sec3.2}. Thus, the results for $\mu_{i}$ under Type-II censoring sample can be obtained using the theorems proved in Subsection \ref{sec3.2}. $(c)$ When $\sigma_{i}$'s are unequal and unknown the joint complete and sufficient statistic is  $(X_{1(1)}, \dots X_{k(1)}, T_1,\dots T_k)$.  Hence, the improved estimators for $\mu_{i}$ will be obtained  similar to Section \ref{sec3.3}.
\end{itemize}
 \subsection{Progressive Type-II censoring}
Suppose in  life-testing  studies,  $n_i$ independents units with lifetimes $X_{i1}, \dots, X_{in_i}$ are drawn from exponential  $Exp(\mu_{i},\sigma_{i})$ distribution. But only $m_i$, $(m_i \le n_i)$ are observed continuously until they fail. Here, censoring takes place progressively in $m_i$ stages. These stages provide failure times of $m_i$ independent units. When the first failure occurs (at the first stage), $S_{i1}$ out of the $ n_i-1$ remaining units are randomly removed from the experiment. Thereafter, at the second failure (at the second stage), $S_{i2}$ out of $ n_i-2-S_{i1}$ remaining units are removed randomly. Eventually when the $m_{i}$ failure occurs (at $m_{i}$ stage) all remaining units $S_{im_{i}} =n_i-m_i-S_{i1}-S_{i2}- \cdots -S_{im_i-1}$ are removed. In this procedure, we have observed the lifetimes of $m_i$ experimental units and the lifetimes are denoted as $X_{i,1:n_i}, \dots X_{i,m_i:n_i}$. These ordered failure times denote the progressive Type-II censored sample. For more information on this censoring scheme, see \cite{balakrishnan2000progressive} and \cite{balakrishnan2007progressive}. Here, our aim is to estimate $\mu_{i}$ based on the progressive Type-II censored samples. From \cite{patra2025inadmissibility}, it can be seen that $X_{i,1:n_i} \sim E(\mu_i,\frac{\sigma_{i}}{n_i})$ and $T_i \sim \Gamma(m_i-1,\sigma_i)$, where $T_i = \sum_{j=i}^{m_{i}}(S_{ij}+1)(X_{i,j:n_{i}}-X_{i,1:n_{i}})$.
\begin{itemize}
\item [(i)] We consider the problem of estimating  $\mu_{i}$ under ordered scale parameter that is $\sigma_{1}  \le \dots \le \sigma_{k}$. Here, $(X_{1,1:n_1},  \dots X_{k, 1:n_k}, T_1,  \dots,  T_k) $ is a joint complete sufficient statistic. Therefore, the present setup reduces to model  (\ref{m1}). Hence, improved estimators for $\mu_{i}$ can be derived similar to that in Section \ref{sec2}.
\item [(ii)]  Now, we consider  estimation of  $\mu_{i}$ under order restriction $\mu_{1} \le \dots \le \mu_{k}$. $(a)$ When $\sigma_{i}$'s are known,  $(X_{1,1:n_1},  \dots, X_{k,1:n_k}) $ is a joint complete and sufficient statistic. This setup aligns with Section \ref{sec3.1}, and consequently improvement results can be obtained similarly. $(b)$  For the case when $\sigma_{i}$'s are equal but unknown, we have joint sufficient statistic as $(X_{1,:n_1}, \dots, X_{k,:n_k}, T)$, which is analogous to Section \ref{sec3.2}. Therefore, improvement results for $\mu_{i}$ under progressive Type-II censoring sample can be proved using the theory derived in Section \ref{sec3.2}. $(c)$ When $\sigma_{i}$'s are unequal and unknown, the joint complete and sufficient statistic is $(X_{1,1:n_1}, \dots X_{k,1:n_k}, T_1,\dots, T_k)$. Thus, the methodologies to obtain improved estimators for $\mu_{i}$ when $\sigma_{i}$'s are unknown under progressive Type-II scheme are similar to Section \ref{sec3.3}.
\end{itemize}
\subsection{Record values}
Record values have wide-range of applications in meteorology, hydrology, and financial market analysis. For generating record values, consider a sequence of independent and identically distributed random variables $Y_1, \dots, Y_{n}$. We define a new sequence of random variables recursively as  
$
Z(1) = 1, \quad \text{and for } k \geq 2, \quad 
Z(k) = \min \{ j \mid j > Z(k-1), \, Y_j > Y_{Z(k-1)} \}.
$
The resulting sequence $\{U_k = Y_{Z(k)}, \, k \geq 1\}$ is known as the sequence of maximal record statistics. For further discussion, please refer to \cite{madi2008improved} and \cite{bobotas2011improved}.
Consider the record sample $X_{i1}^r, \dots, X_{in_{i}}^r $ from exponential $Exp(\mu_{i}, \sigma_i)$ distribution, $i=1,\dots, k$.  From \cite{bajpai2025improved} it has been noted that $
T_{i} = \left( X_{i,n_{i}}^{r} - X_{i,(1)}^{r}\right)  \sim  \Gamma(n_{i}-1, \sigma_{i}), 
$ and $X_{i(1)}^r  \sim Exp(\mu_{i},\sigma_{i})$. Here, we are interested in estimating $\mu_{i}$ under record values for various sceneries as follows.
\begin{itemize}
	\item [(i)] We consider the estimation of $\mu_{i}$ where $\sigma_1\le \dots \le \sigma_k$. In this case, $(X_{i1}^r,\cdots, X_{in_{i}}^r, T_1, \dots  T_k) $ is a joint complete sufficient statistic. Thus, clearly, this model is similar to (\ref{m1}). So improvement results can be obtained similar to Section \ref{sec2}.
\item[(ii)] Now we address the estimation of $\mu_i$ based on the record values when $\mu_1\le \dots \le \mu_k$. Here, we have the following cases.  $(a)$ When $\sigma_{i}$'s are known,  we have $(X_{i1}^r,\cdots, X_{in_{i}}^r) $ is a complete and sufficient statistic. This model is similar to the model discussed in Section \ref{sec3.1}, and consequently we can derive the results for finding the improved estimators for record values using the methodology presented in Section \ref{sec3.1}.  $(b)$  For the case $\sigma_{i}$'s are equal but unknown, we have sufficient statistic as $( X_{i1}^r,\cdots, X_{in_{i}}^r, T)$ with $T=\sum_{i=1}^k T_i$.  This setup is similar to the setup given in Section \ref{sec3.2}. So, the results for finding improved estimators of $\mu_{i}$ can be obtained accordingly as derived in Section \ref{sec3.2}.
$(c)$ When $\sigma_{i}$'s are unequal and unknown, the complete and sufficient statistics for this case are $(X_{i1}^r,\cdots, X_{in_{i}}^r, T_1,\dots T_k)$. This setup is similar to that studied in Subsection \ref{sec3.3}. Therefore, the improved estimators fo  $\mu_{i}$ will be derived using the theories given  in Section \ref{sec3.3}.
\end{itemize}

\section{Conclusions} \label{sec5}
In this work, we have considered component-wise estimation of the location parameter of $k (\ge2)$ exponential distributions under the Linex loss function. We have considered the following scenarios: (i) Estimation of location parameters when scale parameters are ordered, (ii) Estimation of ordered location parameters with  (a) known scale parameter,  (b) equal and unknown scale parameter, and (c) unknown and unequal scale parameters.   For each case, general inadmissibility results have been proved. Using the inadmissibility result for each case, we have obtained  estimator that dominates the usual estimator. Additionally, we have obtained the estimators that dominate MLE and restricted MLE from these results. A detailed numerical comparison of the risk performances of the proposed estimators has been carried out for each case. We have observed that the proposed estimators perform better than the usual estimators. Insights gained from the simulation study have been discussed in detail. Furthermore, it has been noticed that as the value of $p$  increases, the proposed estimators show significant improvement. We have also discussed that using the obtained results one can derive the  results for  Type-II censoring scheme, progressive Type-II censoring scheme and record values. So, This work provides a unified treatment for different important sampling schemes.
\section*{Funding}
Lakshmi Kanta Patra thanks the Science and Engineering Research Board, India for providing financial support to carry out this research  with project number MTR/2023/000229.
\bibliography{bib_oexp}
\end{document}